\documentclass[11pt]{amsart}
\usepackage{palatino}
\usepackage[all]{xy,xypic}
\usepackage{amsfonts,amsmath,oldgerm,amssymb,amscd,tikz-cd,mathrsfs,enumerate}
\usepackage[a4paper,left=12em,right=12em,top=10em,bottom=10em]{geometry}

\makeatletter
\@namedef{subjclassname@2020}{%
	\textup{2020} Mathematics Subject Classification}
\makeatother

\newcommand{\ra}{\rightarrow}		

\newcommand{\by}[1]{\stackrel{#1}{\ra}}

\newcommand{\remove}[1]{}

\newcommand{\surj}{\twoheadrightarrow}	
\newcommand{\inj}{\hookrightarrow}
\newcommand{\ol}{\overline}		
\newcommand{\wt}{\widetilde}
\newcommand{\iso}{\by \sim}

\newtheorem{theorem}{Theorem}[section]
\newtheorem{proposition}[theorem]{Proposition}
\newtheorem{lemma}[theorem]{Lemma}
\newtheorem{definition}[theorem]{Definition}
\newtheorem{corollary}[theorem]{Corollary}
\newtheorem{conjecture}[theorem]{Conjecture}
\newtheorem{example}[theorem]{Example}

\newcommand{\BC}{\mbox{$\mathbb C$}}

	\newcommand{\BN}{\mbox{$\mathbb N$}}

\newcommand{\BS}{\mbox{$\mathbb S$}}

\newcommand{\BR}{\mbox{$\mathbb R$}}

\newcommand{\CI}{\mbox{$\mathcal I$}}	\newcommand{\CJ}{\mbox{$\mathcal J$}}

\newcommand{\CO}{\mbox{$\mathcal O$}}	
	\newcommand{\CR}{\mbox{$\mathcal R$}}

	\newcommand{\CZ}{\mbox{$\mathcal Z$}}

\newcommand{\MI}{\mbox{$\mathfrak I$}}	\newcommand{\MJ}{\mbox{$\mathfrak J$}}
	
\newcommand{\MM}{\mbox{$\mathfrak M$}}

\newcommand{\ma}{\mbox{$\mathfrak a$}}

\newcommand{\mm}{\mbox{$\mathfrak m$}}

\newcommand{\Spec}{\text{Spec}}	
\newcommand{\mS}{\text{MaxSpec}}	
\newcommand{\hh}{\text{ht}}

\newcommand{\rank}{\text{rank}}

\newcommand{\Aut}{\mbox{\rm Aut\,}} 
\newcommand{\Hom}{\text{Hom}\,}

\newcommand{\Um}{\text{Um}}		
\newcommand{\SL}{\text{SL}}
\newcommand{\E}{\text{E}}	
\newcommand{\GL}{\text{GL}}

\newcommand{\bp}{\begin{proposition}}
\newcommand{\ep}{\end{proposition}}
\newcommand{\bl}{\begin{lemma}}
\newcommand{\el}{\end{lemma}}
\newcommand{\bt}{\begin{theorem}}
\newcommand{\et}{\end{theorem}}
\newcommand{\bc}{\begin{corollary}}
\newcommand{\ec}{\end{corollary}}
\newcommand{\bd}{\begin{definition}}
\newcommand{\ed}{\end{definition}}
\newcommand{\bco}{\begin{conjecture}}
\newcommand{\eco}{\end{conjecture}}
\newcommand{\bma}{\begin{bmatrix}}
\newcommand{\ema}{\end{bmatrix}}

\newcommand{\SK}{\mathrm{SK_1}}

\newcommand{\Sp}{\mathrm{K}_1\mathrm{Sp}}

\def\rmk{\refstepcounter{theorem}\paragraph{{\bf Remark} \thetheorem}}
\def\proof{\paragraph{Proof}}
\def\example{\refstepcounter{theorem}\paragraph{{\bf Example} \thetheorem}}
\def\quest{\refstepcounter{theorem}\paragraph{{\bf Question} \thetheorem}}

\title[Zero cycles, Mennicke symbols and $\mathrm{K}_1$-stability of certain real affine algebras]{Zero cycles, Mennicke symbols and $\mathrm{K}_1$-stability\\ of certain real affine algebras}
\author{Sourjya Banerjee}
\address{(Sourjya Banerjee) The Institute of Mathematical Sciences, HBNI, C.I.T. Campus, Tharamani, Chennai  600113, India}
\email{sourjyab@imsc.res.in, sourjya91@gmail.com}
\usepackage{hyperref}

\begin{document}
	
	\maketitle
\begin{abstract}
Let $R$ be a reduced real affine algebra of (Krull) dimension $d \ge 2$ such that either $R$ has no real maximal ideals, or the intersection of all real maximal ideals in $R$ has height at least one. In this article, we prove the following: (1) the $d$-th Euler class group $\E^d(R)$, defined by Bhatwadekar-R.~Sridharan, is canonically isomorphic to the Levine-Weibel Chow group of zero cycles $\text{CH}_0(\Spec(R))$; (2) the universal Mennicke symbol $\text{MS}_{d+1}(R)$ is canonically isomorphic to the universal weak Mennicke symbol $\text{WMS}_{d+1}(R)$; and (3) additionally, if $R$ is a regular domain, then the Whitehead group $\mathrm{SK_1}(R)$ is canonically isomorphic to $\frac{\SL_{d+1}(R)}{\E_{d+1}(R)}$. As an application, we investigate some Eisenbud-Evans type theorems.
\end{abstract}

\section{Introduction}
This article investigates a question of M. P. Murthy on zero cycles \cite[2.12]{M}, the canonical map between the top-length universal weak Mennicke symbol and universal Mennicke symbol, and the stabilization problem for the groups $\mathrm{K}_1$ and $\mathrm{K}_1\mathrm{Sp}$ over certain (mostly singular) real affine algebras. The real affine algebras that are primarily focused on in this article, denoted by the symbol $R$ (throughout the article), satisfy one of the following two conditions, which we call \textbf{P}.
\begin{enumerate}[\quad \quad 1.]
	\item It has no real maximal ideals, or
	\item the intersection of all real maximal ideals has height at least $1$.
\end{enumerate}
Many examples of such real affine algebras are motivated from algebraic geometry. Consider a real affine algebra $A$ such that the closure of the set of all $\mathbb R$-rational points in $\Spec(A)$ has dimension less than $  \dim(\Spec(A))$, then $A$ satisfies \textbf{P}-2. For a concrete example, consider the variety $\Spec(A)$, where $A= \frac{\mathbb R[X_0,\ldots, X_d]}{\langle \sum_{i=0}^{d}X_i^2\rangle}$. Then the only real maximal ideal of $A$ is the maximal ideal corresponding to the origin. This class of real affine algebras is considered in \cite{OPS}, \cite{PO}, and \cite{M}. The main results of this article have three distinct themes, deeply interconnected by the divisibility property of some lower $\mathrm{K}$-groups of certain curves. We discuss each of them separately below.
   
\subsection{Zero cycles and projective modules}
In geometry, the question of determining the precise obstruction for a vector bundle $\mathcal{E}$ on a (connected) affine variety $X$, which can determine whether $\mathcal{E}$ admits a nowhere vanishing section, still remains open. Whenever $\rank(\mathcal{E})>\dim(X)$, J-P. Serre proved in \cite{JPSe} that such a nowhere vanishing section always exist. However, literature provides examples in which J-P. Serre's theorem fails when considering the case $\rank(\mathcal{E})\le \dim(X)$. The only case in which the theory is well developed is the top case where $\rank(\mathcal{E})=\dim(X)$.

In simple algebraic terms, for a commutative noetherian ring $A$ of dimension $d$, the study of a finitely generated projective $A$-module $P$ of rank $d$, along with an obstruction class that can detect the splitting behavior for $P$, remains absent to date, without any further assumption. Nevertheless, a tremendous amount of research has been conducted on this problem since the late 1980s. In \cite{M}, M. P. Murthy settled this question for a smooth affine variety over algebraically closed fields by showing that the top Chern class $\text{c}_d(P)$ in $\text{CH}_0(\Spec(A))$ governs the splitting behavior for $P$. Whenever the affine variety $\Spec(A)$ is not smooth, by $\text{CH}_0(\Spec(A))$ we shall always mean the Chow group of zero-cycles modulo rational equivalence, as defined by M. Levine and C. Weibel in \cite{LW}. It is worth noting that when $\Spec(A)$ is smooth, this definition coincides with the classical definition of the Chow group, e.g. as defined in \cite{Fulton}. Recently, it has been established in \cite{AKMU} that the hypothesis `smoothness' in \cite{M} can be removed. In \cite{JLCTCS}, J.-L. Colliot-Th{\'{e}}l{\`{e}}ne and C. Scheiderer studied the Chow groups of real affine varieties. Although, over an arbitrary base field, the vanishing of the top Chern class is a necessary condition for $P$ to split off a free factor, it is not sufficient, even for smooth real varieties. The study of a sufficient condition on a smooth real affine variety for the splitting problem of top rank projective modules via its top Chern class was initiated in \cite{BRS}. It took a deep understanding of two other groups, namely the Euler class group $\E^d(A)$ and the weak Euler class group $\E^d_0(A)$ defined in \cite{SMBB3} to completely solve the problem in \cite{BMDSM}. A diligent reader might notice that the following theorem plays a crucial step in this development when dealing with smooth real varieties.

\bt\cite[Theorem 5.5]{BRS}
Let $\Spec (A)$ be a smooth affine variety of dimension $d\ge 2$
over $\BR$. Then the canonical surjective map $\E_0^d(A)\iso \text{CH}_0(\Spec(A))$ is an isomorphism.
\et
Therefore, the following question naturally arises when dealing with singular real affine algebras. Later in this subsection, we briefly mention how it plays a crucial role in understanding many open questions in literature.
  
\smallskip

\quest Let $A$ be a reduced real affine algebra of dimension $d \geq 2$. Is the canonical surjective map $\E_0^d(A) \surj \text{CH}_0(\Spec(A))$ an isomorphism?
 
 \smallskip
 
The above question remains open, except in the smooth case. Recently, A. Asok and J. Fasel proved in \cite{AF} that the aforementioned question has an affirmative answer for a smooth affine scheme over an infinite perfect field of characteristic $\neq 2$. In Theorem \ref{QBRS} we prove the following. 

\bt
The canonical surjective map $\E_0^d(R) \iso \text{CH}_0(\Spec(R))$ is an isomorphism, where $\dim(R)=d\ge 2$.
\et

We give a sketch of our approach. The first step of the proof is to show that the group $\E^d(R)$ is divisible [Proposition \ref{decg}]. To prove this, we utilize another variant of Euler class groups $\widetilde{\E}^d(R)$ and $\widetilde{\E}^d_0(R)$ [see \ref{DECG}], which are closely related to \cite[Definition 2.1]{BRS} and  \cite[7.3]{RGAK}. These two definitions of the Euler class groups coincide when the base ring is a smooth affine algebra over an infinite perfect field. The notation $\mu(-)$ stands for the minimal number of generators. In the second step, we employ the divisibility of $\E^d(R)$ to establish that any reduced ideal $J\subset R$ such that $\hh(J)=\mu(J/J^2)=d$ is a surjective image of a projective $R$-module $P$ of rank $d$. Furthermore, we compute its Chern class with respect to an ideal $K$ supported only on smooth complex maximal ideals [for details we refer to Proposition \ref{PG2}]. In the third step, using results in \cite{JLCTCS} and \cite{AKMU} we show that the cycle $[K]$ associated to $R/K$ in $\text{CH}_0(\Spec(R))$ is not a torsion element. Then applying the second step it follows that the canonical map is an isomorphism. 

As two of the interesting consequences, we (i) prove an analogy to A. A. Rojtman's theorem \cite{Roj}, and (ii) provide an affirmative answer to \cite[Question 2]{MKD1}. In the following, we describe a few of these consequences [for details we refer to \ref{1C}].

\bt
Let $\dim(R)=d\ge 2$, and let $P$ be a projective $R$-module of rank $d$. Let $I\subset R$ be an ideal of height $d$ such that $\mu(I/I^2)=d$. Then the following are true.

\begin{enumerate}[\quad \quad (1)]
	\item $\text{CH}_0(\Spec(R))$ is torsion free.
%	\item  $P$ splits into a free factor if and only if  $C_d(P)=0$ in $\text{CH}_0(\Spec(R))$.
	\item  Suppose that $\ol f:P/IP\surj I/I^2$ be a surjective map. Then there exists a surjective lift  $f:P\surj I$ of $\ol f$ if and only if $e(P)=(I)$ in $\E^d_0(R,L)$, where $\wedge^dP=L$.
	\item Let $J\subset R[T]$ be an ideal such that $\mu(J/J^2)=\hh(J)=d$. Then there exists a projective $R[T]$-module $P$ of rank $d$ such that $J$ is a surjective image of $P$. 
	\item Let $R(T)$ be the ring obtained from $R[T]$ by inverting all monic polynomials in $R[T]$. Then for all $d\ge 3$ the canonical map $\Gamma :\E^d(R[T])\to \E^d(R(T))$ is injective.
\end{enumerate}
\et

\subsection{Mennicke Symbols} The study of the elementary orbit space of unimodular rows $\frac{\Um_{d+1}(A)}{\E_{d+1}(A)}$ is well-documented in the literature, where $A$ is a ring of dimension $d\ge 2$. L. N. Vaser{\v{s}}te{\u{\i}}n, for $d=2$ \cite[Section 5]{SV} and W. van der Kallen, for $d\ge 3$ \cite{VdK}, have shown that the orbit space $\frac{\Um_{d+1}(A)}{\E_{d+1}(A)}$ possesses an abelian group structure. Furthermore, it is proven in \cite{VdKM} that the group $\frac{\Um_{d+1}(A)}{\E_{d+1}(A)}$ coincides with the universal weak Mennicke symbol $\text{WMS}_{d+1}(A)$ of length $d+1$. From the definition, one may observe that there exists a canonical surjection $f:\text{WMS}_{d+1}(A) \twoheadrightarrow \text{MS}_{d+1}(A),$ where $\text{MS}_{d+1}(A)$ is the universal Mennicke symbol of length $d+1$. Since the Bass-Kubota theorem is true for one-dimensional rings \cite[Theorem 2.12]{VdK}, and considering the inductive approach of W. van der Kallen, as used in establishing the group structure in $\frac{\Um_{d+1}(A)}{\E_{d+1}(A)},$ we are naturally led to the following question: Is $f$ an isomorphism? An example over real affine algebras provides evidence that this is not the case in general (cf. \cite[4.17]{VdK} and \cite[Example 2.2(c)]{CAW}).

The group $\text{WMS}_{d+1}(A)$ is said to have a \textit{nice} group structure if for any two unimodular rows $(a,a_1,\ldots,a_d)$ and $(b,a_1,\ldots,a_d)$ of length $d+1$, the following holds. $$[(a,a_1,\ldots,a_d)]\star[(b,a_1,\ldots,a_d)]=[(ab,a_1,\ldots,a_d)]$$ Here $[-]$ denotes elementary orbit space of $-$. One may observe that $f$ is an isomorphism if and only if $\text{WMS}_{d+1}(A)$ has a nice group structure (cf. \cite[Theorem 2.1]{AR}). The injectivity of the map $f$ has been studied across various algebras in a series of articles, including \cite{AR}, \cite{Fasel}, \cite{AAR}, \cite{SK} and \cite{SK3}. In Theorem \ref{NGS}, we prove that the group $\text{WMS}_{d+1}(R)$ has a nice group structure, where $\dim(R)=d\ge 2$. As a consequence, in Corollary \ref{NGSC} we show that $\text{WMS}_{d+1}(R)$ is a divisible group.

\subsection{Stability for $\mathrm{K}_1$ and $\mathrm{K}_1\mathrm{Sp}$ groups} We begin this part by recalling the stabilization problem for the Whitehead group $\mathrm{SK_1}(A)$, where $A$ is a ring. It follows from the definition of $\mathrm{SK_1}$ that there exists a canonical map $\Gamma_n:\frac{\SL_{n}(A)}{\E_n(A)}\to \mathrm{SK_1}(A)$, for all $n\ge 2$. The stabilization problem for the Whitehead group $\mathrm{SK_1}(A)$ asks the following. 
\smallskip

\quest Given a ring $A$ of dimension $d\ge 2$, what is the least positive integer $n$ for which the maps $\Gamma_{n+i}$ become injective (similarly surjective) for all $i\ge 0$?
\smallskip

A similar question can be asked for the symplectic group $\mathrm{K_1Sp}$ with the appropriate adjustments. For an arbitrary commutative noetherian ring, the bounds on such an $n$ beyond which both surjectivity and injectivity hold were established by L. N. Vaser{\v{s}}te{\u{\i}}n in \cite{V69}, \cite{Va1} and \cite{Vas}. One can construct examples of smooth real varieties, which can establish that L. N. Vaser{\v{s}}te{\u{\i}}n's bounds are optimal in general (cf. \cite[Proposition 7.10]{VdKM}, \cite{RN}, \cite{FSR}). However, it was proved in \cite{RV}, \cite{BRAO} and \cite{BCR} that L. N. Vaser{\v{s}}te{\u{\i}}n's injective stability bound for both $\mathrm{K}_1$ and $\mathrm{K}_1\mathrm{Sp}$ can be improved for regular affine algebras over a field $k$ of characteristic zero such that the cohomological dimension of $k$ is at most one (the condition on the base field $k$ is more general, see e.g. \cite[Proposition 3.1]{RV}). We prove analogues of their theorems in the following form [for details we refer to Section \ref{3}].

\bt
In addition to condition \textbf{P}, if we assume that $R$ is a regular domain, and let $I\subset R$ be a principal ideal, then the following are true.
\begin{enumerate}[\quad\quad (1)]
	\item $\mathrm{SK_1}(R,I)\cong\frac{\SL_{n}(R,I)}{\E_{n}(R,I)}$, for all $n\ge d+1$.
	\item Suppose that $d\not\equiv 2\mod 4$. Then the canonical map $\frac{\mathrm{Sp}_{2n}(R,I)}{\mathrm{ESp}_{2n}(R,I)}\to \mathrm{K_1Sp}(R,I)$ is injective for all $n\ge [\frac{d+1}{2}]$.
\end{enumerate}
 
\et 

\subsection{Some Eisenbud-Evans type theorems} Section \ref{6} is an extension of Section \ref{1}. Here we investigate relations between some Eisenbud-Evans type theorems as studied in \cite{NMK84}. Let $A$ be a ring. Recall that the group $\mathrm{F^dK}_0(A)$, which is the subgroup of $\mathrm{K_0}(A)$ generated by the images of all $[A/I]$, where $I$ is a locally complete intersection ideal of height $d$. One of the main results in Section \ref{6} is Theorem \ref{k_0}, where we prove that $\mathrm{F^dK}_0(R)$ is canonically isomorphic with $\E^d(R)$. This provides an affirmative answer to M. P. Murthy's torsion problem \cite[2.12]{M} for $R$.

\subsection{Notations}\label{SR1} Unless otherwise stated, all rings considered in this article are assumed to be commutative noetherian with $1(\not=0$), and all modules are assumed to be finitely generated. Any projective module is assumed to have a constant rank. We denote the vector $(1, \ldots, 0)$ as $e_1$. For a ring $A$, the set $\mS(A)$ is the collection of all maximal ideals in $A$. Throughout the article the symbol $R$ will always represent a real affine algebra of dimension $d\ge 2$ satisfying one of the following two conditions \textbf{P}.

\begin{enumerate}[\quad \quad 1.]
	\item It has no real maximal ideals, or
	\item the intersection of all real maximal ideals has height at least $1$.
\end{enumerate}

\section*{Acknowledgment}This work was initiated when I was a PhD student at the Indian Statistical Institute, Kolkata, under the supervision of Mrinal Kanti Das. Some parts of Sections \ref{2} and \ref{3} were included in my PhD thesis. I am grateful to Mrinal Kanti Das for introducing me to this class of real algebras and for his invaluable guidance throughout my research journey. I also thank Kuntal Chakraborty for helping in locating a couple of references on symplectic matrices and Mennicke symbols. I am grateful to Rahul Gupta for assisting me with a couple of references included in the revised version. I am deeply indebted to the referees for careful reading of the manuscript and for pointing out a gap in Corollary \ref{pfpb} in an earlier version. Without their detailed comments, the exposition would have lacked clarity at various places.

\section{Preliminaries}
This section is devoted to recollecting several results and definitions from the literature that serve as prerequisites for the main theorems in this article. In some cases, we slightly modify or restate well-known results to better suit our requirements.

\subsection{Stably free modules}
The purpose of this subsection is to give the proof of Theorem \ref{P}. We believe Theorem \ref{P} must be well known due to A. A. Suslin \cite{AAS} and others (see \cite{PO}). However, we could not find an exact reference. We begin with the following definition.

\bd
Let $A$ be a ring and $M$ be an $A$-module. The \textit{order ideal} of an element $m\in M$ is defined by $\CO_M(m):=\{\alpha(m)\in A:\alpha\in M^*=\Hom_R(M,A)\}.$ An element $m\in M$ is called a \textit{unimodular element} if $\CO_M(m)=A$. The set of all unimodular elements of $M$ is denoted by $\Um(M)$. If $M=A^n$, then we write $\Um_n(A)$ instead of $\Um(A^n)$.
\ed

We recall Swan's version of Bertini theorem from \cite[Theorem 1.3]{SBE}.

\bt\label{SwB} Let $V=\Spec(A)$ be a smooth affine variety over an
infinite field $k$. Let $Q$ be a projective $A$-module of rank $r$ and $(q, a)\in \Um(Q \oplus A)$. Then there is a $y\in Q$ such that the ideal
$I: = \CO_Q(q + ay)$ has the following properties.
\begin{enumerate}[\quad\quad (1)]
	\item The subscheme $U=\Spec (A/I)$ of $V$ is smooth over $k$ and $\dim (U)=
	\dim (V)-r$, unless $U=\phi$;
	\item If $\dim (U) \not= 0$ then $U$ is a variety.
\end{enumerate}
\et

We now recall the following definitions from \cite[Chapter 1]{SV}.

	\bd Let $A$ be a ring and $I$ be an ideal of $A$. We define the following.
\begin{enumerate}
	\item The \textit{elementary group} $\E_{n}(A)$ is the subgroup of $\GL_{n}(A)$ generated by the matrices $\E_{ij}(x) := \text{I}_n + \text{e}_{ij}(x)$, where $x\in A$, $i \neq j$, and $\text{e}_{ij}(x)$ is the matrix with only possible non-zero entry $x$ at the $(i,j)$-th position.
	\item $\Um_{n}(A,I):=\{v\in \Um_{n}(A):v\equiv e_1\mod I\}$, i.e., $\Um_{n}(A, I)$ the set of all \textit{relative unimodular rows} of length $n$ with respect to the ideal $I$.
	
	\item $\SL_{n}(A,I):=\{\alpha\in \SL_{n}(A):\alpha\equiv \text{I}_n\mod I \}$. 
	
	\item For all $n\ge 2$, the group $\E_n(A,I)$ is the smallest normal subgroup
	of $\E_n(A)$ containing the elements $\E_{21}(x)$, where $x \in I$ (see \cite[2.1]{VdK}).
%	\item Two matrices $\alpha\in \text{SL}_r(A,I)$ and $\beta\in  \text{SL}_s(A,I)$ are said to be \textit{K-equal}, if 
%	
%	$$\begin{pmatrix}
%		\alpha & 0\\
%		0 & \text{I}_s
%	\end{pmatrix}=\begin{pmatrix}
%		\beta & 0\\
%		0 & \text{I}_r
%	\end{pmatrix}.$$The relation K-equal is an equivalence relation.
	\item For any two matrices $\alpha\in \text{M}_r(A)$ and $\beta\in  \text{M}_s(A)$ we define $$\alpha\perp \beta:=\begin{pmatrix}
	\alpha & 0\\
	0 & \beta
\end{pmatrix} \in \text{M}_{r+s}(A).$$

	\item One can identify $\SL_{n}(A,I)$ as a subgroup of $\SL_{n+1}(A,I)$ by associating a matrix $\alpha\in \SL_n(A,I)$ with the matrix $1\perp \alpha \in \SL_{n+1}(A,I)$. We define the set $\text{SL}(A,I):=\bigcup_n\text{SL}_n(A,I) $. Note that, the set $\text{SL}(A,I)$ forms a group in an obvious way, which coincide with the group structure on $\text{SL}_r(A,I)$ for all $r\ge 1$. 
	\item We define $\text{E}(A,I):=\bigcup_n\text{E}_n(A,I)$ as a subgroup of $\text{SL}(A,I) $. 
	\item We define $\mathrm{SK_1}(A,I):=\frac{\text{SL}(A,I)}{\E(A,I)}$.

	\item We inductively define an alternating matrix $\chi_r\in \E_{2r}(A)$ as follows:
	$$\chi_1:=\begin{pmatrix}
		0 & 1\\
		-1 & 0
	\end{pmatrix}\in \E_2(A)\text{ and }\chi_{r+1}:=\chi_r\perp \chi_1.$$
	\item We define $\text{Sp}_{2n}(A):=\{\alpha\in \SL_{2n}:\alpha^T\chi_n\alpha=\chi_n\}$.
	\item We define $\text{Sp}_{2n}(A,I):=\{\alpha\in \text{Sp}_{2n}(A):\alpha\equiv \text{I}_{2n}\mod I\}$.
	%		\item Let $\sigma$ be the permutation of the natural numbers given by $\sigma(2i) = 2i-1$ 
	%		and $\sigma (2i-1) = 2i$ for $i = 1, 2,..., n$.

	\item Let $\sigma$ denote the permutation of the natural numbers given by $\sigma(2i)=2i-1$
	and $\sigma(2i-1) = 2i$. Let $x\in A$. Then $\text{sE}_{ij} (x)$ is defined as follows:
	$\text{sE}_{ij}(x):=\text{I}_{2n} + \text{e}_{ij}(x)$ if $i = \sigma(j)$; and 
	$\text{sE}_{ij} (x):=\text{I}_{2n} + \text{e}_{ij}(x) - (-1)^{i+j} \text{e}_{\sigma(j)\sigma(i)}(x)$ if $i \not= \sigma(j)$ and $i<j$.
	\item The subgroup of $\text{Sp}_{2n}(A)$ generated by $\text{sE}_{ij} (x)$ is called the elementary symplectic group $\text{ESp}_{2n}(A)$, where $x\in A$ and $i\not =j$.
	\item The group $\text{ESp}_{2n}
	(A, I)$ is defined to be the smallest normal
	subgroup of $\text{ESp}_{2n}
	(A)$ containing $\text{sE}_{ij} (x)$, where $x\in I$ and $i\not =j$.
	\item We define $\text{Sp}(A,I):= \bigcup_n\text{Sp}_{2n}(A,I)$ as a subgroup of $\text{SL}(A,I) $.
	\item We define $\text{ESp}(A,I):= \bigcup_n\text{ESp}_{2n}(A,I)$ as a subgroup of $\text{SL}(A,I) $.
	\item We define $\mathrm{K_1Sp}(A,I):=\frac{\text{Sp}(A,I)}{\text{ESp}(A,I)}.$
	\item A matrix $\alpha\in \GL_{n}(A)$ is called \textit{isotopic to identity} if there exists $\theta(T)\in \GL_{n}(A[T])$ such that $\theta(0)=\text{I}_n$ and $\theta(1)=\alpha$.
\end{enumerate}
\ed

The next two results follow from a theorem due to A. A. Suslin \cite{AAS} and their proofs can be found in \cite[Propositions 3 and 4]{PO}.

\bp\label{po}
Let $C$ be a smooth real curve having no real maximal ideal. Then $\SK(C)$ is a divisible group.
\ep

\bp\label{pa}
Let $C$ be as in Theorem \ref{po}. Then the natural map $\Sp(C)\iso \SK(C)$ is an isomorphism.
\ep

	\smallskip
	The proof of the following lemma can be found in \cite[Corollary 2.3]{AAS}.

	\bl\label{AASC}
	
	Let $A$ be a ring and $n\geq 2$. Suppose that $(a_0,\ldots,a_n) \in \Um_{n+1}(A)$ satisfies $\dim(A/\langle a_2,\ldots,a_n\rangle )\leq 1$ and $\dim(A/\langle a_3,\ldots,a_n\rangle )\leq 2$. Furthermore, assume that there exists $\alpha \in \SL_2(A/\langle a_2,\ldots,a_n\rangle )\cap \mathrm{ESp}(A/\langle a_2,\ldots,a_n\rangle )$ such that $(\overline{a_0},\overline{a_1})\alpha =(\overline{b_0},\overline{b_1})$. Then, there exists $\gamma\in \E_{n+1}(A)$ such that $(a_0,\ldots,a_n)\gamma=(b_0,b_1,a_2,\ldots,a_n)$.
	\el
	
	The next proposition is based on a clever observation of R. Parimala, that one may avoid singularities on A. A. Suslin's proof of \cite[Theorem 2.4]{AAS}. This proposition plays a crucial role in this article. Hence, we give a sketch of the proof. 
	
	\bp\label{PR}
	Let $A$ be a reduced affine algebra of dimension $n$, over a perfect field $k$, and let $v=(v_0,\ldots,v_n)\in \Um_{n+1}(A)$. Assume that $\BS\subset \mS(A)$ is such that the ideal $ \CI:=\bigcap\limits_{m\in\BS} m$, has height $\ge 1$. Then there exists a matrix $\epsilon\in \E_{n+1}(A)$ such that if we take $((u_0,\ldots, u_n)=)u=v\epsilon$ then
	
	\begin{enumerate}[\quad \quad (1)]
		\item $A/\langle u_0,\ldots,u_{i-1}\rangle $ is a smooth affine $k$-algebra of dimension $n-i$, and
		\item  $\mS(A/\langle u_0,\ldots,u_{i-1}\rangle )\cap \BS=\emptyset$, for all $1\le i\le n$.
%		\item $A/\langle u_0,\ldots,u_{n-1}\rangle $ is a product of fields.
	\end{enumerate}
	% $(i)$ $u_{0}-1\in I$, where $I$ is the ideal defining the singular locus of $A.$
	Additionally, we may assume that $A/\langle u_0,\ldots,u_{i-1}\rangle $ is a smooth affine domain if $i<n$.
	
	\ep
	
	\proof  Since $A$ is reduced, and $k$ is a perfect field, the algebra $A$ is geometrically reduced by \cite[\href{https://stacks.math.columbia.edu/tag/020I}{Lemma 33.6.3}]{TSP}. Therefore, it follows from \cite[\href{https://stacks.math.columbia.edu/tag/056V}{Lemma 33.25.7}]{TSP} that the regular locus of $A$ contains all generic points of $\Spec(A)$. Hence, the ideal $\CJ$ defining the singular locus of $A$ has height at least $ 1$. Let $I= \CI \CJ$, then by our hypothesis we obtain that $\hh(I)\ge 1$. Hence, going modulo $I$ and using standard stability arguments (such as the Prime avoidance lemma), one can find an $\widetilde{\omega}\in \E_{n+1}(A/I)$ such that $v\widetilde{\omega} \equiv e_1\mod I$. We choose a lift $\Omega\in \E_{n+1}(A)$ of $\widetilde{\omega}$. Then we have $v\Omega = w=(w_0,\ldots,w_{n})$, where $1-w_0\in I$ and $w_i\in I$ for all $i\ge 1$. Since it is enough to prove the theorem for $w$, we may begin with the assumption that $1-v_0\in I $ and $v_i\in I$ for $i\ge 1$. 
	
	Also, we observe that it is enough to prove the theorem for $i=1$, as then we can repeat the same inductive steps on each of the rings $A/\langle v_0,v_1,\ldots, v_{t}\rangle $ for $t=1,2,\ldots, i-2$, to obtain the result on $A/\langle v_0,\ldots,v_{i-1}\rangle $. Moreover, at each step, since the completion of the corresponding unimodular row is elementary, we can always lift it back to our original ring $A$. Furthermore, the fact that $1-v_0\in I\subset \CI $ will ensure that $\BS\cap \mS(A/\langle v_0,\ldots , v_{i-1}\rangle )=\emptyset$. Hence, in the remaining part, we only prove the case $i=1$. 
	
	Applying Theorem \ref{SwB}, we get $\lambda_j\in A$, for $j=1,\ldots,n$, such that if we replace $v_0$ by $v'_0=v_0+\sum\lambda_jv_j$, then $\Spec (A/\langle v'_0\rangle )$ is a smooth irreducible variety outside the singularities of $A$ such that $\dim(A/\langle v'_0\rangle )=n-1$. We take $u_0=v_0'$ and $u_i=v_i$ for $i\ge 1$. Then there exists an $\epsilon\in \E_{n+1}(A)$ such that $u=v\epsilon$. We now show that $u$ satisfies (1) and (2). As $1-v_0\in I$ and $v_i\in I$ for $i\ge 1$, we have $u_0-1\in I$. Since $A/\langle u_0\rangle $ is a domain (in particular reduced), the set of all singular points in $\Spec (A/\langle u_0\rangle )$ forms a closed set. Let $J$ be the ideal defining the singular locus of $A/\langle u_0\rangle $. Since $A/\langle u_0\rangle$ is smooth outside the singularities of $A$ we have $\ol I\subset J$, where `bar' denotes going modulo $\langle u_0\rangle$. However $u_0-1\in I$ gives us the fact that $\ol I=A/\langle u_0\rangle$. Hence $A/\langle u_0\rangle$ is smooth and  $\mS(A/\langle u_0\rangle )\cap \BS=\emptyset$. This concludes the proof.\qed

	Now we are ready to prove the main theorem of the subsection.
	
	\bt\label{P} 
	Any unimodular row in $R$ of length $d+1$ is completable to the first row of an invertible matrix. As a consequence, all stably free $R$-modules of rank $d$ are free.
	\et 
	\proof First we note that if $d<2$, then there is nothing to prove. Hence, we may assume that $d\ge 2$. One may observe that, if $R$ satisfies \textbf{P}-2, then this is done in \cite[Theorem 3.2]{OPS}. Therefore, we assume that $R$ is a real affine algebra having no real maximal ideal. Let $v=(v_0,\ldots,v_d)\in \Um_{d+1}(R)$. Using Lemma \ref{PR}, we may assume that $C:=R/\langle v_0,\ldots,v_{d-2}\rangle $ is a smooth curve. Let `bar' denote going modulo $\langle v_0,\ldots,v_{d-2}\rangle $. 

	Since by Theorem \ref{po} the group $\SK(C)$ is divisible, there exists a $\sigma\in \SL_2(C)\cap \E_3(C)$ such that $\sigma(\ol v_0,\ol v_1)=(\ol a^{d!},\ol b)$. Applying Theorem $\ref{pa}$, we get $\sigma\in \SL_2(C)\cap \mathrm{ESp}(C)$. Now one may use Lemma \ref{AASC} to find an $\epsilon\in \E_{d+1}(R)$ such that $v\epsilon=(a^{d!}, b,v_3,\ldots,v_d)$. It follows from \cite{AASUS} that the unimodular row $(a^{d!}, b,v_3,\ldots,v_d)$ is completable to the first row of an invertible matrix.\qed

%%%%%%%%%%%%%%%%%%%%%%%%% Alternative proof using Suslin's %%%%%%%%%%%%%%%%%%%%%%%%%%%%%%%%%%%%
%	First we note that if $d<2$, then there is nothing to prove. Hence we may assume that $d\ge 2$. Note that, if $R$ satisfies P-2, then this is done in \cite[Theorem 3.2]{OPS}. Therefore, it is enough to assume that $R$ is a real affine algebra having no real maximal ideal. Let $v=(v_0,\cdots,v_d)\in \Um_{d+1}(R)$. Using Lemma \ref{PR}, we may assume that $C:=R/\langle v_0,\cdots,v_{d-2}\rangle $ is a smooth curve. Let `bar' denote going modulo $\langle v_0,\cdots,v_{d-2}\rangle $. 
%	
%	Since by Theorem \ref{po} the group $\SK(C)$ is divisible, there exists a $\sigma\in \SL_2(C)\cap \E_3(C)$ such that $\sigma(\ol v_0,\ol v_1)=(\ol a^{d!},\ol b)$. Applying Theorem $\ref{pa}$, we get $\sigma\in \SL_2(C)\cap \mathrm{ESp}(C)$. Now one may use Lemma \ref{AASC} to find an $\epsilon\in \E_{d+1}(R)$ such that $v\epsilon=(a^{d!}, b,v_3,\cdots,v_d)$. It follows from \cite{AASUS} that the unimodular row $(a^{d!}, b,v_3,\cdots,v_d)$ is completable to the first row of an invertible matrix.\qed 
%%%%%%%%%%%%%%%%%%%%%%%%%%%%%%%%%%%%%%%%%%%%%%%%%%%%%%%%%%%%%%%%%%%%%%%%%%%%%%%%%%%%%%%%%%%%%%%

	\rmk \label{fr} For any unimodular row $w=(w_0,\ldots,w_{d})\in \Um_{r+1}(A)$, we shall call $(w_0,\ldots,w_{r-1}, w_r^{r!})$ as the factorial row of $w$ and denote it by $\chi_{r!}(w)$. It follows from the proof of Theorem \ref{P} that, for any $v\in \Um_{d+1}(R)$ there exist an $\epsilon\in \E_{d+1}(R)$ and $w\in \Um_{d+1}(R)$ such that $v\epsilon=\chi_{d!}(w)$.

	\subsection{Cancellation of projective modules}
	Let $A$ be a ring. Recall that a projective $A$-module $P$ of rank $n$ is said to be \textit{cancellative} if $P \oplus A^r \cong P' \oplus A^r$ implies $P \cong P'$, for any  $A$-module $P'$. After Theorem \ref{P}, a natural question arises: whether a projective $R$-module of rank $d$ is cancellative. It turns out that this question has an affirmative answer. Again we believe this result is well-known. However, we were unable to find a suitable reference for the exact following version. When $R$ satisfies condition \textbf{P}-2, then the result follows from \cite[Theorem 3.2]{OPS}. We sketch a proof of the next theorem, which is along similar lines to the proof of Theorem \ref{P}. Before that, we recall the following definition.
	
	\bd
	Let $A$ be a ring. Let $P$ be a projective $A$-module such that $\Um(P)\not=\phi$. We choose $\phi\in P^*$ and $p \in P$ such that $\phi(p)=0$. We define an endomorphism $\phi_p$ as the composite $\phi_p:P\to A\to P$, where $A\to P$ is the map sending $1\to p.$ Then by a transvection we mean an automorphism of $P$, of the form $1+\phi_p$, where either $\phi\in \Um(P^*)$ or $p\in \Um(P)$. By $\E(P)$  we denote the subgroup of $\Aut(P)$ generated by all transvections.
	\ed
	\rmk If $P$ is a free $A$-module of rank $n\ge 3$, then it is proved in \cite[Lemma 2.20]{BBR} that $\E(P)$ coincides with $\E_n(A)$.
	
	\bt\label{PC}
	Let $P$ be a projective $R$-module of rank $d$. Then $P$ is cancellative.
	\et
	\proof  It is sufficient to assume that $R$ satisfies \textbf{P}-1. Furthermore, without loss of generality, we may assume that $R$ is reduced. Let us choose $(a,p)\in \Um(R\oplus P)$. We shall show that there exists $\sigma\in \Aut(A\oplus P)$ such that $\sigma(a,p) = (1,0)$. Let $J$ be the ideal defining the singular locus of $R$, then $\hh(J) \geq 1$. There exists a non-zero divisor $t\in J$ such that $P_t$ is free. We choose a set of generators of $P_t$. Then by clearing denominators, we can find an element $s = t^l$ for some $l \in \mathbb{N}$ such that $sP \subset R^d \subset P$.
	
	Since $s$ is a non-zero divisor, applying a classical result due to Bass (\cite{HBASS} or see \cite[Proposition 2.13]{SMB}), we obtain that the canonical map $\Um(R\oplus P)\surj \Um(\frac{R}{\langle s \rangle}\oplus \frac{P}{ s P})$ is surjective. Therefore, without loss of generality, we may assume that $a-1\in\langle s\rangle$ and $p\in sP\subset R^d$. Hence, we may take $p=(a_1,\ldots,a_d)\in R^d$.
	
	Using Lemma \ref{PR}, we may further assume that $B=A/\langle a,a_1,\ldots, a_{d-2}\rangle $ is a smooth affine domain of dimension one. Since $P_t$ is free and $a-1\in \langle t\rangle$, we obtain that the module $ \frac{P}{\langle a,a_1,\ldots,a_{d-2}\rangle P}$ is a free module over $B$. Since $P/aP$ is a free module over $R/\langle a \rangle$, we have $\E_d(\frac{R}{\langle a\rangle} )= \E(\frac{P}{\langle a \rangle P})$. Now we may follow the arguments given in the proof of Theorem \ref{P} to obtain a $\gamma \in \E_d(R/\langle a\rangle )$ such that $\gamma (\widetilde{a_1},\ldots,\widetilde{a_d})=(\widetilde{a_1},\ldots,\widetilde{a_{d-2}},\widetilde{b_{d-1}},\widetilde{b_d}^{d!})$, where `tilde' denotes going modulo $\langle a\rangle $. Since $\langle a,s\rangle =R$ and $sP\subset R^d\subset P$, we obtain that $\widetilde P=(\frac{R}{\langle a \rangle})^d$. In particular, we get the equality $\Um(\widetilde{P})=\Um_{d}(R/\langle a\rangle)$. Applying \cite[Proposition 2.12]{SMB} we can lift $\gamma\in \E_d(R/\langle a\rangle )$ to an $\alpha\in \Aut(P)$. Therefore, we obtain that $$\alpha p\equiv({a}_1,\ldots,{a}_{d-2},{b}_{d-1},{b}_d^{d!})\mod \langle a\rangle   P.$$ We define $\sigma :=\begin{pmatrix}
		1 & 0\\
		0 & \alpha
	\end{pmatrix}$. Moreover, as $p\in sP\subset R^d\subset P$, we get the following: $$ (a,p)\equiv (a,\alpha p)\equiv (a,a_1,\ldots,a_{d-2},b_{d-1},b_d^{d!})\mod \E(R\oplus P).$$ Since factorial rows are completable by \cite{AASUS}, one may conclude the proof (cf. \cite[Theorem 4.1]{SMB}).\qed
	
	\subsection{Excision ring and relative cases} The purpose of this section is to state and prove a relative version of Theorem \ref{P}, which will be needed to improve the injective stability of the Whitehead group $\SK$ in Section \ref{3}. We begin by recalling some interesting facts about the excision ring (for the definition, we refer to \cite[Definition 2.5]{AAR}). Let $A$ be a ring. For any ideal $I\subset A$, the excision ring $A\oplus I$ can be viewed as a fiber product of $A$ with respect to the ideal $I$. Therefore, if $A$ is an affine algebra over a field $k$, then $A\oplus I$ is also an affine algebra over $k$ \cite[Proposition 3.1]{MKA}. We denote $\pi_2:A\oplus I\surj A$ by the map sending $(r,i)\mapsto r+i$. The following lemma ensures that the excision ring $R\oplus I$ enjoys similar properties as $R$.
	
	\bl\label{EL}
	Let $I\subset R$ be an ideal. Then $R\oplus I$ is also a real affine algebra of dimension $d$ satisfying condition \textbf{P}.
	\el
	\proof Firstly, we note that $R\oplus I$ is also a real affine algebra of dimension $d$, and the extension $R\inj R\oplus I$ is an integral extension, for a proof we refer to \cite[Proposition 3.1]{MKA}. Therefore, it suffices to show that $R\oplus I$ satisfies condition \textbf{P}. We show this in separate cases.
	
	Let $R$ satisfy \textbf{P}-1, and let $\MM\subset R\oplus I$ be a maximal ideal. Then $\frac{R}{R\cap \MM}\inj \frac{R\oplus I}{\MM}$ is an integral extension. Therefore, the ideal $R\cap \MM$ is a maximal ideal in $R$. Since $R$ has no real maximal ideal, we have $\frac{R}{R\cap \MM}\cong \mathbb C$. Hence $\frac{R\oplus I}{\MM}$ is also isomorphic to $\mathbb C$.
	
	Let $R$ satisfy \textbf{P}-2, and let $\MJ$ be the intersection of all real maximal ideals in $R\oplus I$. Let $\MM$ be a real maximal ideal in $R\oplus I$. Since $R\inj R\oplus I$ is an integral extension, we have $\frac{R}{R\cap \MM}\cong \mathbb R$. This fact combined with the fact that the intersection of all real maximal ideals of $R$ has height $\ge 1$ implies that $\hh(\MJ\cap R)\ge 1$. Let $a\in \MJ\cap A$ be a non-zero divisor. Then the element $(a,0)\in \MJ$ is a non-zero divisor in $R\oplus I$. This concludes the proof. \qed

	%First, we note that any ideal of $R\oplus I$ is of the form $J\oplus I'$, where $J\subset R$ and $I'\subset I$ are ideals of $R$. Thus, in particular, maximal ideals of $R\oplus I$ are of the form $m\oplus I$, where $m\subset R$ is a maximal ideal. Now if $R$ satisfies condition P-1 then the residue field of $R\oplus I$ remains $\mathbb C.$ So if $R$ satisfy P-1, then so does $R\oplus I$.
	%
	%Now suppose assume $R$ is satisfying P-2. Let $\CJ$ be the intersection of all real maximal ideals then $\hh(\CJ)\ge 1$. Then we note that the intersection of all real maximal ideals of $R\oplus I$ is $\CJ\oplus I$ which also has a positive height. Hence $R\oplus I$ satisfies condition P-2.\qed

Now we are ready to state and prove the main result of this subsection. The proof essentially uses the idea of \cite[Proposition 3.3]{RV}.
%\bt
%Let $\sigma\in \SL_{d+1}(R)$ be stably elementary matrix, then $\sigma$ is isotopic to identity.  Furthermore assume $R$ to be smooth. Then $E_{d+2}(R)\cap \SL_{d+1}(R)=E_{d+1}(R)$, for $d\ge 2.$
%\et
%\proof Let $\sigma\in E_{d+2}(R)\cap \SL_{d+1}(R)$. Get $\tau(T)\in E_{d+2}(R[T])$ such that $\tau(0)=I_{d+2} $ and $\tau(1)=1\perp\sigma.$ Let $t=T^2-T\in R[T]$ be a non zero element and $v=e_1\tau(T)\in \Um_{d+2}(R[T])$. Thus $v\cong e_1$ modulo $\langle t\rangle $ in $R[T]$. Thus by Theorem \ref{SS1} there exists $\chi(T)\in \SL_{d+2}(R[T])$, such that $\chi(T)\cong I_{d+2}$ modulo $\langle t\rangle $ and $v=e_1\chi(T)$.\\
%Thus $e_1\tau(T)\chi(T)^{-1}=e_1 $. Hence $\tau(T)\chi(T)^{-1}$ is of the form $(1\perp \rho(T))\prod_{i=1}^{d+2}E_{i,1}(\lambda_i)$, where $\lambda_i\in R[T]$ and $\rho(T)\in \SL_{d+1}(R[T])$.\\
%Now since $\chi(T)\cong I_{d+2}$ modulo $\langle t\rangle $, we have $\chi(0)=\chi(1)=I_{d+2}$. Thus $\rho(0)=I_{d+1}$ and $\rho(1)=\sigma$ that is $\rho$ is an isotopy of $\sigma.$\\
%Now let $d \ge 1$. By  Vorst’s $K_1$-analogue in \cite{Bo}
%of H. Lindel’s theorem in \cite{Li}, $\SL_{d+1}(R_p[T]) = E_{d+1}(R_p[T])$ for any prime
%ideal $p \in \Spec(R)$. By the “Local–Global Principal” for $E_r(R[T])$, $r \ge 3$
%in \cite{Ra}, $\rho(T) \in \rho(0)E_{d+1}(R[T])$, and so $\rho(T) \in E_{d+1}(R[T])$ as $\rho(0) = I_{d+1}$.
%But then $\sigma = \rho(1) \in E_{d+1}(R)$. \qed

\bp\label{S1}
Let $v\in \Um_{d+1}(R,I)$, where $I=\langle a\rangle$ for some non-zero divisor $a\in R$. Then $v$ can be completed to a matrix $\sigma\in \SL_{d+1}(R,I)$.
\ep

\proof Let $v=(1+av_1,av_2,\ldots, av_{d+1})$. Since $v\equiv e_1 \mod \langle a \rangle$, there exists $u\in \Um_{d+1}(R)$ such that $vu^t=1$ and $u\equiv e_1\mod\langle a \rangle$. Let $u=(1+au_1,au_2,\ldots,au_{d+1})$. Since $a$ is a non-zero divisor we obtain that $u_1+v_1=-a\sum_{i=1}^{d+1}u_iv_i$. Let $B=R[T]/\langle T^2-aT\rangle $, and let `bar' denote going modulo $\langle T^2-aT\rangle$. We take $v(T)=(1+v_1T,v_2T,\ldots, v_{d+1}T)$ and $u(T)=(1+u_1T,u_2T,\ldots, u_{d+1}T)$. Then we get the following. $$v(T)u(T)^t=1+T(u_1+ v_1)+ T^2{(\sum_{i=1}^{d+1}u_iv_i)}= 1+T(T- a){(\sum_{i=1}^{d+1}u_iv_i)}\equiv 1 \mod\langle T^2-aT\rangle$$ This shows that $\ol{v(T)}\in \Um_{d+1}(B)$. We observe that $B\cong R\oplus \langle a\rangle$ \cite[Corollary 3.2]{MKA}. Hence, by Lemma \ref{EL} it follows that $B$ also satisfies condition \textbf{P}. Therefore, by Theorem \ref{P}, there exists an $\alpha( T)\in \SL_{d+1}(B)$, such that \begin{equation}\tag{$\star$}  \ol{v(T)}=e_1\alpha(T). \end{equation}Using the identifications $\frac{B}{\langle\ol T\rangle }\cong R$ and $\frac{B}{\langle \ol T-\ol a\rangle} \cong R$, we obtain that $\ol{v(0)}=e_1$ and $\ol{v(a)}=v$. Hence from ($\star$) we get $e_1=e_1\alpha(0)$ and $v=e_1\alpha(a)$. Let $\sigma=\alpha(0)^{-1}\alpha(a)$. Then $e_1\sigma=v$, and $\sigma\equiv \text{I}_{d+1}$ mod $\langle a\rangle $. This concludes the proof. \qed

The following lemma can be proven using exactly the same arguments given in \cite[Theorem 4.4]{AG}. One just needs to use Remark \ref{fr} in place of \cite[Lemma 4.1]{AG} in their proof. Hence, we skip the proof.

\bl\label{bla}
Let $a\in R$ be a non-zero divisor and let $I=\langle a \rangle$. Let $v\in \Um_{d+1}(R,I)$. Then there exist $w\in \Um_{d+1}(R,I)$ and $\epsilon\in \E_{d+1}(R,I)$ such that $v\epsilon=\chi_{d!}(w)$.
\el

%	\bt\label{SS1}
%	Let $d\ge 2$ and let $I\subset R$ be an ideal. Then $\Um_{d+1}(R,I)=e_1\SL_{d+1}(R,I).$
%	\et
%	\proof Let $v=(v_0,\cdots,v_d)\in \Um_{d+1}(R,I)$. Then we note that $$\tilde{v}:=((1,v_0-1),(0,v_1)\cdots,(0,v_d))\in \Um_{d+1}(R\oplus I,0\oplus I).$$Applying Theorem \ref{P} and Lemma \ref{EL} we get $\alpha\in \SL_{d+1}(R\oplus I)$ such that $e_1\alpha=\tilde{v}$. Let `bar' denote going modulo $0\oplus I$. Then we have $\ol e_1\ol \alpha=\ol e_1$, where $\ol \alpha \in \SL_{d+1}(R)\subset \SL_{d+1}(R\oplus I)$. Replacing $\alpha$ by $ \ol \alpha^{-1}\alpha $ one may further assume that $\alpha\in \SL_{d+1}(R\oplus I, 0\oplus I).$ Here we note that this replacement does not change the fact that $e_1\alpha =\tilde{v}$. Then $e_1\SL_{d+1}(\pi_2)(\alpha)=v$, where $\SL_{d+1}(\pi_2):\SL_{d+1}(R\oplus I)\to \SL_{d+1}(R)$ induced by $\pi_2$. Now note that we have the following commutative diagram
%	
%	$$\begin{tikzcd}
%		1 \arrow{rr} && \SL_{d+1}(R\oplus I,0\oplus I)\arrow{rr}\arrow[d,dashrightarrow] && \SL_{d+1}(R\oplus I) \arrow{rr}\arrow{d}{\SL_{d+1}(\pi_2)} && \SL_{d+1}(\frac{R\oplus I}{0\oplus I})\arrow{d}{\ol{\SL_{d+1}(\pi_2)}} \\
%		1 \arrow{rr} && \SL_{d+1}(R, I)\arrow{rr} && \SL_{d+1}(R) \arrow{rr} && \SL_{d+1}(\frac{R}{I}) 
%	\end{tikzcd}$$ 
%	Thus $\SL_{d+1}(\pi_2)$ and $\ol{\SL_{d+1}(\pi_2)}$ will induce $\Gamma: \SL_{d+1}(R\oplus I,0\oplus I)\to \SL_{d+1}(R,I)$ such that the diagram commutes. Hence we actually get $\SL_{d+1}(\pi_2)(\alpha)\in \SL_{d+1}(R,I)$ and this completes the proof. \qed
	
	\subsection{Mennicke symbols}
 In this subsection, we briefly recall the Mennicke symbols and the weak Mennicke symbols. We begin with the following definitions. 
	\bd
	Let $A$ be a ring. A Mennicke symbol of length $n + 1 \ge 3$, is a pair $(\psi, G)$,
	where $G$ is a group and $\psi : \Um_{n+1}(A) \to G$ is a map such that:
	\begin{enumerate}[\quad\quad (1)]
		\item  $\psi((0,...,0, 1)) = 1$ and $\psi(v) = \psi(v\epsilon)$ for any $\epsilon \in \E_{n+1}(A)$;
		\item $\psi((b_1,..., b_n, x))\psi((b_1,...,b_n, y)) = \psi((b_1,...,b_n, xy))$ for any two unimodular rows $(b_1,..., b_n, x)$ and $(b_1,...,b_n, y)$.
	\end{enumerate}
	\ed
	
	\bd
	Let $A$ be a ring. A weak Mennicke symbol of length $n +1\ge 3$ is a
	pair $(\psi,G)$, where $G$ is a group and
	$\phi : \Um_{n+1}(A) \to G$
	is a map such that:
	\begin{enumerate}[\quad \quad (1)]
		\item $ \phi(1, 0, ..., 0) = 1$ and $\phi(v) = \phi(w\epsilon)$ if $\epsilon \in \E_{n+1}(A)$;
		\item $\phi(a, a_1,..., a_n) \phi(1-a, a_1,... , a_n) = \phi(a(1-a), a_1, ..., a_n)$ for any unimodular row $(a, a_1,..., a_n)$ such that $(1-a, a_1,..., a_n)$ is also unimodular.
	\end{enumerate}
	\ed
	
		Let $A$ be a ring of dimension $n$. Clearly, a universal Mennicke symbol $(\text{ms}, \text{MS}_{n+1}(A))$ and a universal weak Mennicke symbol $ (\text{wms},\text{WMS}_{n+1}(A))$ exist. With an abuse of notation, we only use the notations $\text{MS}_{n+1}(A)$ and $\text{WMS}_{n+1}(A)$. In \cite{VdK}, W. van der Kallen defined an abelian group structure on the elementary orbit space of unimodular rows $\frac{\Um_{n+1}(A)}{\E_{n+1}(A)}$. Moreover, in \cite{VdKM} it was shown that the group $\frac{\Um_{n+1}(A)}{\E_{n+1}(A)}$ coincides with the group $ \text{WMS}_{n+1}(A)$. Hence, we will stick to the notation $ \text{WMS}_{n+1}(A)$ only.

	\subsection{Euler class groups}\label{DECG}

	 In this subsection, we give a slightly alternative definition of top Euler class groups for $R$. Our definition is motivated from \cite[Definition 2.1]{BRS} and  \cite[7.3]{RGAK}. In Section \ref{1} we show that our definition matches with the Euler class groups defined in \cite{SMBB3} for $R$.
	 
	In addition to condition \textbf{P}, let $R$ be reduced. Let $I\subset R$ be an ideal such that $\mu(I/I^2)=\hh(I)=d$. Let $L$ be a rank one projective $R$-module. Let $\alpha,\beta:(L/IL)\oplus(R/I)^{d-1}\surj I/I^2$ be two surjections. We say $\alpha$ and $\beta$ are related if there
	 exists $\sigma\in \E((L/IL)\oplus (R/I)^{d-1})$ such that $\alpha \sigma = \beta$. This defines an equivalence relation on the set of all surjections $(L/IL)\oplus(R/I)^{d-1}\surj I/I^2$. Let $[\alpha]$ denote the
	 equivalence class of $\alpha$. We shall call $[\alpha]$ a \textit{local $L$-orientation} of $I$. With abuse of notation sometimes we shall call $\alpha$ a local $L$-orientation of $I$ instead of $[\alpha]$. The local $L$-orientation $[\alpha]$  is said to be a \textit{global $L$-orientation} of $I$, if there exists a surjective map $\Gamma:L\oplus R^{n-1}\surj I$ such that $\Gamma\otimes R/I=\alpha$. In this situation, we call $\Gamma$ a surjective lift of $\alpha$. Notice that, since the canonical map $\E(L\oplus R^{d-1})\surj \E((L/IL)\oplus (R/I)^{d-1})$ is surjective, if $\alpha$ has a surjective lift, then so has any $\beta$ equivalent to $\alpha$. Whenever $L\cong R$, we shall say a local orientation (respectively global orientation) of $I$ instead of an $R$-local orientation (respectively $R$-global orientation) of $I$.
	 
	 %		Let $\dim(R)=d\ge 2$. Let $I\subset R$ be an ideal such that $\mu(I/I^2)=\hh(I)=d$.

	 Let $\widetilde{G}$ be the free abelian group on the set $B$ of pairs $(\mm, 
	 \omega_{\mm} )$, where:
	 \begin{enumerate}
	 	\item $\mm \subset R$ is a smooth complex maximal ideal of height $d$;
%	 	\item $\mm$ is a $\MM$-primary ideal, for some maximal ideal $\MM\subset A$ of height $d$;
	 	\item $\omega_{\mm} : (L/\mm L)\oplus (R/\mm)^{d-1}\surj \mm/\mm^2$
	 	is a local $L$-orientation of $\mm$.
	 \end{enumerate}
	 
	  We shall say that an ideal $J \subset R$ of height $d$ is \textit{regular} if it
	 is of the form $J = \cap_i \mm_i$, where each $\mm_i$ is a smooth complex maximal ideal of height $d$. Given a regular ideal $J=\bigcap_i \mm_i$, we observe that there always exists a local $L$-orientation. Let $\omega_J$ be a local $L$-orientation of $J$. Then $\omega_J$ gives rise, in a natural way, to a local $L$-orientation $\omega_{\mm_i}$ of $\mm_i$. We associate to the pair
	 $(J, \omega_J)$, to the element $\sum_i (\mm_i,\omega_{\mm_i})$ of $\widetilde{G}$. By abuse of notation, we denote the element $\sum_i (\mm_i,\omega_{\mm_i})$ by $(J,\omega_J)$.
	 
	 Let $\widetilde{H}$ be the subgroup of $\widetilde{G}$ generated by the set $S$ of pairs $(J, \omega_J )$, where $\omega_J$ is a global $L$-orientation of $J$. Then the quotient group $\widetilde{G}/\widetilde{H}$ is the \textit{$d$-th Euler class group} of $R$ with respect to $L$, denoted as $\widetilde{\E}^d(R,L)$. Whenever $L\cong R$, we shall write $\widetilde{\E}^d(R)$ instead of $\widetilde{\E}^d(R,R)$.

%	 Let $P$ be a projective $A$-module of rank $d$ with determinant $L$. Let $\chi:\wedge^d(L\oplus A^{d-1})\iso \wedge^d (P)$ be an isomorphism. Applying Swan's version of Bertini theorem \cite[Theorem 1.3]{SBE} one can find a surjection $\lambda : P \surj I$ such that $I$ is a product of finitely many distinct maximal ideals in $A$ of height $d$. Let `bar' denote going modulo $I$. We obtain an induced surjection $\lambda\otimes (A/I) : P/IP \surj  I/I^2$. Note that, since $P$ has determinant $L$ and
%	 $\dim (A/I) =0$, we have $P/IP\cong (L/IL)\oplus (A/I)^{d-1}$. We choose an isomorphism $\phi :(L/IL)\oplus (A/I)^{d-1}\iso  P/IP$, such that $\wedge^d\phi=\chi\otimes A/I$. Let $\omega_I$ be the surjection $ (\lambda\otimes A/I)\circ\phi: (L/IL)\oplus(A/I)^{d-1}\surj I/I^2$. We say that $(I, \omega_I)$ is an ``Euler cycle'' induced by the triplet
%	 $(P,\lambda,\chi)$. 
%	 
%	 Whenever the Euler cycle is independent of the choice of $\lambda$, that is, if $(J,\omega_J)$ is another Euler cycle induced by a triplet $(P,\lambda',\chi)$, for some $\lambda'\in P^*$, then  $(I,\omega_I)=(J,\omega_J)$ in $\widetilde{\E}^d(R,L)$, in this situation we shall call the class of $(I,\omega_I)$ in the Euler class group $\widetilde{\E}^d(A,L)$ as the ``Euler class" of the pair $(P,\chi)$, denoted as $e(P,\chi)$.
	 
%	  It was proved in \cite{SMBB3}, that whenever the ring contains an infinite field and either $A$ is smooth or $\frac{1}{d!}\in A$, then the Euler class of such a pair $(P,\chi)$ is well defined.
	 
	 Similarly, one can define the $d$-th weak Euler class group $\widetilde{\E}^d_0(R,L)$ as $\widetilde{F}/\widetilde{K}$, where $\widetilde{F}$ and $\widetilde{K}$ are as follows. 
	 
	 Let $\widetilde{F}$ be the free abelian group on the set of all smooth complex maximal ideals $\mm$ of height $d$. Let $J=\bigcap_i \mm_i$ be a regular ideal. Then we associate $(J)$ to the element $\sum_i (\mm_i)$ of $\widetilde{F}$. By abuse of notation, we denote the element $\sum_i (\mm_i)$ by $(J)$.
	 
	 Let $\widetilde{K}$ be the subgroup of $\widetilde{F}$ generated by $(J)$ such that $J$ has a global $L$-orientation. Whenever $L\cong R$, we shall write $\widetilde{\E}_0^d(R)$ instead of $\widetilde{\E}_0^d(R,R)$. 
	 
	We denote the Euler class groups and weak Euler class groups defined in \cite{SMBB3} by $\E^d(R,L)$ and $\E_0^d(R,L)$ respectively. One may observe that from the universal property of quotient group there exist canonical homomorphisms $i_R:\widetilde{\E}^d(R,L)\to \E^d(R,L)$ and $\phi_R : \widetilde{\E}^d_0(R,L) \to \text{E}^d_0(R,L)$.

	\subsection{Chow groups} Let $A$ be a reduced affine algebra of dimension $d\ge 1$, over a field and $X=\Spec(A)$. In this subsection, we briefly recall the definition of the Chow group $\mathrm{CH}_0(X)$ as defined by M. Levine and C. Weibel in \cite[$\S$1]{LW}. We would like to mention that their original definition is more general compared to the one given below. Here, we restrict it to our specific set-up.

	Let $X_{\text{reg}}$ be the disjoint union of smooth loci of the $d$-dimensional irreducible components of $X$. Let $X_\text{sing}$ denote the complement of $X_\text{reg}$ in $X$. Let $\CZ_0 (X)$ be the free abelian group on closed points lying in $X_\text{reg}$.

	A \textit{Cartier curve} of $X$ is a closed affine subscheme $Z\subset X$ of dimension $1$ such that (1) no component of $Z$ is contained in $X_\text{sing}$, and (2) every point of $Z\cap X_\text{sing}$ lies in a neighborhood $U$ such that $Z\cap U$ is defined
	by a regular sequence of functions on $U$.

Let $C$ be a Cartier curve of $X$. Let $\CZ_0 (C,C\cap X_\text{sing})$ be the free abelian group on all closed points lying in $C \cap X_\text{reg}$. We shall often write  $\CZ_0 (C,C\cap X_\text{sing})$ as  $\CZ_0 (C)$. Then there is a natural inclusion $\CZ_0 (C)\inj \CZ_0 (X)$.

 A section $f$ of $\CO_C$ is said to be a \textit{rational function} on $C$ if $f$ belongs to the total quotient ring $k(C)= \prod_c \CO_{C,c}$, where the product runs for all generic points $c$ of $C$. The group $k(C)^{*}$ is the subgroup of units of $k(C)$ which are units in $\CO_{C,x}$, for every $x\in C\cap  X_{\text{sing}}$. Recall that, there exists a well-defined cycle map $\mathrm{cyc}_C:k(C)^*\to \CZ_0(C)$, as defined in \cite[Definition 1]{LW}. Then $\CR_0(X)$ is the subgroup of $\CZ_0(X)$ generated by the images of the following compositions of maps
 
 $$\begin{tikzcd}
k(C)^* \arrow{r}{\mathrm{cyc}_C}  & \CZ_0(C) \arrow[hookrightarrow]{r} & \CZ_0(X),
\end{tikzcd}$$where $C$ runs over all Cartier curves of $X$.  Then the \textit{Chow group} is defined as the quotient $$\text{CH}_0(X):= \frac{\CZ_0(X)}{\CR_0(X)}.$$

	In addition to condition \textbf{P}, let $R$ be reduced. Let $\mm\in \mS(R)$ be a regular maximal ideal. Then the assignment $\mm \mapsto [\mm] \in  \text{CH}_0(\Spec(R))$ induces a canonical
	group homomorphism $\widetilde{F} \to  \text{CH}_0(\Spec(R))$, where the cycle $[\mm]$ is associated to the class of $R/\mm$ in $\text{CH}_0(\Spec(R))$. It follows from \cite[Lemma 2.2]{LW} that this assignment kills the classes of complete intersection ideals. Therefore, it defines a natural group homomorphism $\widetilde{\E}_0^d(R) \to \text{CH}_0(\Spec(R))$. 
	
%	\bd Let $A$ be a ring. Let $I\subset A$ be an ideal such that $\mu(I/I^2)=n$. The ideal $I$ is said to be projective generated if there exists a projective $A$-module $P$ of rank $n$ such that $I$ is a surjective image of $P$.
%	\ed 

\section{Zero cycles and projective modules}\label{1}
	In addition to condition \textbf{P}, let $R$ be reduced. Let $P$ be a projective $R$-module of top rank. One of the primary objectives of this section is to show that the $d$-th Euler class group $\E^d(R)$ is canonically isomorphic with the Chow group $\text{CH}_0(\Spec(R))$. As a preparation, we prove various results that are interesting in their own right. We begin with a rephrased version of a moving lemma within our framework. This arises as a direct consequence of the Bertini-Murthy-Swan theorem \cite[Corollary 2.6]{M}. We shall employ this lemma repeatedly throughout this section. Hence, we give a sketch.

\bl\label{SB}
In addition to condition \textbf{P}, let $R$ be reduced. Let $J\subset R$ be an ideal containing a non-zero divisor such that $\mu(J/J^2)=d$. Then there exists an ideal $I\subset R$ such that 
\begin{enumerate}[\quad \quad (1)]
	\item $I$ is a product of finitely many distinct smooth complex maximal ideals of height $d$,
	\item $I+JK=R$, for any ideal $K$ such that $\hh(K)\ge 1$, and
	\item $\mu(I\cap J)=d$.
\end{enumerate}
Moreover, if $\hh(J)=d$, then one can choose $I$ with the following additional properties. 
\begin{enumerate}[\quad \quad (4)]
	\item Let $L$ be a projective $R$-module of rank one. Then, for any given local $L$-orientation $\omega_J$ of $J$, we get a local $L$-orientation $\omega_I$ of $I$ and a local $L$-orientation $\omega_{IJ}$ of $IJ$ such that $(IJ,\omega_{IJ})=(I,\omega_{I})+(J,\omega_{J})=0$ in $\E^d(R,L)$. Consequently, we get $(IJ)=(I)+(J)=0$ in $\E_0^d(R,L)$. 
\end{enumerate}
\el

\proof Let $J=\langle f_1,\ldots, f_d\rangle +J^2$, and $K$ be an ideal in $R$ such that $\hh(K)\ge 1$. Let $J_1$ be the ideal defining the singular locus of $R$. Let $J_2$ be the intersection of all real maximal ideals in $R$ if $R$ satisfies \textbf{P}-2, else we take $J_2=R$. We define $\MI:=J_1\cap J_2\cap K$, then $\hh(\MI)\ge 1$. Now we apply  Bertini-Murthy-Swan theorem \cite[Corollary 2.6]{M} to get $h_i\in J$ and an ideal $I$, which is a product of distinct smooth maximal ideals of height $d$ such that $IJ=\langle h_i,\ldots, h_d\rangle$, $I+J\MI=R$, and $f_i-h_i\in J^2$. Moreover, since $I+J_2=R$ the ideal $I$ is only supported on smooth complex maximal ideals. Hence this proves (1), (2) and (3).

Now let us assume that $\hh(J)=d$. The local $L$-orientation $\omega_J$ induces a map ${f}:L\oplus R^{d-1}\to J$, such that $f(L\oplus R^d)+J^2=J$. Now, we repeat the arguments given in the first paragraph and obtain an ideal $I$ which is a product of finitely many distinct smooth complex maximal ideals of height $d$ and a surjection $g:L\oplus R^{d-1}\surj  IJ$ with $f\otimes R/J \equiv g\otimes R/J$ such that $I+\MI J=R$. Let $\omega_{IJ}$ be the global $L$-orientation of $IJ$ induce by $g$. Now since $I+J=R$, by Chinese remainder theorem $g$ will induce local $L$-orientations $\omega_I$ of $I$ and $\omega_J'$ of $J$. Moreover, as $f\otimes R/J \equiv g\otimes R/J$, we obtain that $(J,\omega_J)=(J,\omega_J')$ in $\E^d(R,L)$. This in particular gives us that $(IJ,\omega_{IJ})=(I,\omega_{I})+(J,\omega_{J})=0$ in $\E^d(A,L)$. Now the remaining part of (4) follows from using the canonical group homomorphisms $\E^d(R,L)\to \E_0^d(R,L)$. This concludes the proof.\qed

We now proceed to prove the divisibility of the Euler class group $\E^d(R)$. To that end, we begin with a series of preparatory results, starting by recalling the following computation from \cite[Lemma 4.1]{DM2}.

\bl\label{BDML}
Let $I,J\subset R$ be two ideals of height $d$, and let $n\in \mathbb N$. Suppose that, there exist $a_i\in R$ such that $I=\langle a_1,\ldots,a_d\rangle +I^2$ and $J=\langle a_1,\ldots,a_{d-1}\rangle +I^{n}$. Let $L$ be a projective $R$-module of rank one. Then $(J)=n(I)$ in $\E_0^d(R,L).$
\el

\bl \label{Step1}
Let $L$ be a projective $R$-module of rank one. The canonical map $i_R:\widetilde{\E}^d(R,L)\iso\E^d(R,L)$ is an isomorphism.
\el

\proof We note that if $(I,\omega_I)\in \widetilde{\E}^d(R,L)$ such that its image vanishes in $\E^d(R,L)$, then by \cite[Theorem 4.2]{SMBB3} it follows that $\omega_I$ is a global $L$-orientation of $I$, and hence $(I,\omega_I)=0$ in $\widetilde{\E}^d(R,L)$. This shows that $i_R$ is injective.

Let $(J,\omega_J)\in \E^d(R,L)$, where $J\subset R$ is an ideal such that $\mu(J/J^2)=\hh(J)=d$. Then by Lemma \ref{SB} there exists an ideal $I$ which is a product of finitely many distinct smooth complex maximal ideals of height $d$ and a local $L$-orientation $\omega_I$ of $I$ such that $(J,\omega_J)=-(I,\omega_I)$ in $\E^d(R,L)$. Since $I$ is a product of distinct smooth complex maximal ideals we get $(I,\omega_I)\in \widetilde{\E}^d(R,L)$, and hence $(J,\omega_J)\in \widetilde{\E}^d(R,L)$. This completes the proof. \qed 

\bl\label{Step2}
Let $L$ be a projective $R$-module of rank one. The canonical map $\psi_R:\widetilde\E^d(R,L)\iso \E^d_0(R,L)$ sending $(\mm, 
\omega_{\mm} )\mapsto (\mm)$ is an isomorphism, where $\mm\subset  R$ is a smooth complex maximal ideal  of height $d$.
\el

\proof We observe that, by definition, the canonical map $\tau_R: \E^d(R,L) \to \E^d_0(R,L)$ is surjective. Moreover, it follows again from the definitions that the following triangle is commutative.

$$\begin{tikzcd}
	\E^d(R,L) \arrow{r}{\tau_R}  & \E^d_0(R,L) \\
	\widetilde\E^d(R,L) \arrow{u}{i_R} \arrow{ur}[swap]{\psi_R} & 
\end{tikzcd}$$
Since $i_R$ is an isomorphism by Lemma \ref{Step1}, we get $\psi_R$ is surjective. Hence, we only prove the injective part. 

%	Let $(I,\omega_I)\in \E^d(R,L)$ be such that there exists a global $L$-orientation $\omega_I'$ of $I$. First, we show that $(I,\omega_I)=0$ in $\E^d(R,L)$. Since by Theorem \ref{PC} the projective module $L\oplus R^{d-1}$ is cancellative this part is standard in the literature. Hence we only give a sketch.

Let $(I,\omega_I)\in \widetilde\E^d(R,L)$ be such that there exists a global $L$-orientation $\omega_I'$ of $I$. First, we show that $(I,\omega_I)=0$ in $\widetilde\E^d(R,L)$. We note that two local $L$-orientations of $I$ must differ by a unit in $R/I$. Let $\omega_I$ and $\omega_I'$ differ by $u\in (R/I)^*$. Since $I$ is a product of finitely many distinct smooth complex maximal ideals, there exists $v\in (R/I)^*$ such that $u=v^2$. Since $\omega_I'$ is a global $L$-orientation, applying \cite[Lemma 5.3]{SMBB3} we conclude that $(I,\omega_I)=0$ in $\E^d(R,L)$. Therefore, by \cite[Theorem 4.2]{SMBB3} it follows that $\omega_I$ is a global $L$-orientation of $I$. Hence $(I,\omega_I)=0$ in $\widetilde{\E}^d(R,L)$.

Now let $(J,\omega_J)$ be in the kernel of the map $\widetilde{\E}^d(R,L)\to \E^d_0(R,L)$. Then the same proof of \cite[Lemma 3.3]{BRS} yields the following
$$(J,\omega_J)+\sum_{i=1}^{r}(J_i,\omega_i)=\sum_{k=r+1}^{s}(J_k,\omega_k) \text{ in } \widetilde{\E}^d(R,L),$$ where $J_i,J_k$ are regular ideals in $R$ of height $d$ such that there exist surjections $R^{d-1}\oplus L\surj J_j$, for each $j$. However, then in the previous paragraph we have already proved that $(J_j,\omega_j)=0$, for each $j$. Hence $(J,\omega_J)=0$ in $\widetilde{\E}^d(R,L)$. This concludes the proof.\qed

Combining Lemmas \ref{Step1} and \ref{Step2}, we obtain the following corollary.

\bc\label{Cor}
Let $L$ be a projective $R$-module of rank one.	The canonical group homomorphism $\Phi_R:\text{E}^d(R,L)\iso\text{E}^d_0(R,L)$ is an isomorphism. 
\ec

\bl\label{Step3}
Let $L$ be a projective $R$-module of rank one.	The canonical group homomorphism $\phi_R : \widetilde{\E}^d_0(R,L) \iso \text{E}^d_0(R,L)$ is an isomorphism.
\el

\proof Let $I\subset R$ be a regular ideal of height $d$ such that $\phi_R(I)=0$ in $\text{E}_0^d(R,L)$. It follows from Corollary \ref{Cor} and \cite[Theorem 4.2]{SMBB3} that any local $L$-orientation of $I$ is a global $L$-orientation. Therefore, we get $(I)=0$ in $\widetilde{\E}^d_0(R,L)$, and hence $\phi_R$ is injective.

To prove the surjectivity, we choose an element $(J)\in \E^d_0(R,L)$, where $J\subset R$ is an ideal such that $\mu(J/J^2)=\hh(J)=d$. Now applying Lemma \ref{SB} there exists an ideal $K\subset R$ of height $d$ such that

\begin{enumerate}[\quad \quad (1)]
	\item $\mu(K\cap J)=d$;
	\item $K+J=R$;
	\item $K$ is a product of distinct smooth maximal ideals in $R$ of height $d$.
\end{enumerate}
The conditions (1) and (2) imply that $(J)+(K)=0$ in $\E^d_0(R,L)$. From (3) it follows that $(K)\in \widetilde{\E}^d_0(R,L)$, and hence $\phi_R$ is surjective.\qed

We are now ready to prove the divisibility of the Euler class group $\text{E}^d(R)$.

	\bp\label{decg}
	The Euler class group $\E^d(R)$ is divisible.
	\ep
	\proof Note that, by \cite[Corollary 4.6]{SMBB3} the group $\E^d(R)$ is canonically isomorphic to the group $\E^d(R_{\text{red}})$. Hence, without loss of generality, we may assume that $R$ is reduced. Let $L$ be a projective $R$-module of rank one. It follows from Corollary \ref{Cor} and Lemma \ref{Step3} that it suffices to prove the divisibility of the group $\widetilde{\E}^d_0(R)$. 
	
	 	 Let $m\subset R$ be a smooth complex maximal ideal of height $d$, and let $n\in \mathbb N$. It follows from Lemma \ref{Step3} that it is enough to find an ideal $I\subset R$ such that $\mu(I/I^2)=\hh(I)=d$, and $(m)=n(I)$ in $\E^d_0(R)$. We now follow the argument of M. P. Murthy \cite{M}.
	 	
	 	Let $J_1$ be the ideal defining the closed set which is the
	 	union of the singular locus of $\Spec (R)$ and all the irreducible components of
	 	$\Spec (R)$ of dimension less than d. Let $J_2$ be the intersection of all real maximal ideals in $R$, if $R$ satisfies \textbf{P}-2, else we take $J_2=R$.  We set $J:=J_1\cap J_2$. By the choice of $m$ it follows that $m+J=R$. Let $m=\langle f_1, \ldots, f_d\rangle +m^2$. Then applying \cite[Corollary 2.4]{M} there exist $g_i\in m^2$, $i=1,\ldots, d-1$ such that if we set $h_i=f_i+g_i$, then \begin{enumerate}[\quad \quad (1)]
	 		\item $\langle h_1,\ldots, h_{d-1} \rangle +J=R$, and
	 		\item  $C:=R/\langle h_1,\ldots, h_{d-1} \rangle $ is a regular ring of dimension $1$.
	 	\end{enumerate}  We note that, it follows from (1) that the curve $C$ does not have any real maximal ideals. Hence $\text{Pic}(C)$ is a divisible group. Furthermore, we have the following. $$m=\langle h_1, \ldots,h_{d-1}, f_d\rangle +m^2$$ Let `bar' denote going modulo $\langle h_1,\ldots, h_{d-1}\rangle $. As $\mu(\ol m/\ol m^2)=1$, the ideal $\ol m$ is an invertible ideal in $C$. Since $\text{Pic}(C)$ is a divisible group, there exists an ideal $\ol K\subset C$ such that	$(\ol m)^{-1}=(\ol K)^n  \text{ in }  \text{Pic}(C)$. Moreover, using moving arguments without loss of generality we may assume that $\ol m+\ol K=C$, and $\ol K$ is a product of distinct maximal
	 	ideals in $C$. Then we can find a non-zero divisor $\ol h\in C$ such that the following holds.\begin{equation}\tag{$\star$} \ol m\ol K^n=\langle \ol h\rangle 	\end{equation}Let $K\subset R$ be a lift of $\ol K$. Then $ K$ is a product of distinct smooth complex maximal ideals in $R$ of height $d$. Hence $\mu(K/K^2)=\hh(K)=d$. Moreover, as $\ol K$ is an invertible ideal we have $\mu(\ol K/\ol K^2)=1$. Let $\ol K=\langle \ol h_d \rangle + \ol K^2$, and let $h_d\in R$ be a lift of $\ol h_d$. Then $K=\langle h_1,\ldots, h_{d} \rangle +K^2$. We take $K'=\langle h_1,\ldots, h_{d-1} \rangle + K^n$. Then from Lemma \ref{BDML} we get $(K')=n(K)$ in $\E^d_0(R)$. It follows from $(\star)$ that $\mu(m\cap K')=d$. Hence $(m)=-(K')=-n(K)$ in $\E^d_0(R)$. Therefore, we take $(I)=-(K)$ to conclude the proof.\qed

	\smallskip
	
%	\rmk\label{inv} Applying \cite[Definition 3.2]{YY} the extension $R\inj R_{\BC}$ will induce a map $\chi_R:\E^d(R)\to \E^d(R_{\BC})$. One may observe that the map $\Gamma_R $ is an isomorphism, and the inverse of $\Gamma_R$ is the map induced by $\chi_R$.

	\rmk We would also like to point out an obvious observation of Lemmas \ref{Step1} and \ref{Step2}. Let $P$ be a projective $R$-module of rank $d$ with determinant $L$, and $I\subset R$ be an ideal such that $\mu(I/I^2)=d=\hh(I)$. Applying \cite[Corollary 4.4]{SMBB3} we get that (1) $(I)=0$ in $\E_0^d(R,L)$ if and only if $\mu(I)=d$, and (2) $e(P)=0$ in $\E_0^d(R,L)$ if and only if $P$ has a unimodular element.

%\rmk It follows from [Proposition \ref{decg}, Step - 1, 2]	and \cite[Theorem 6.8]{SMBB3} that the groups $\E^d(R)$ and $\E^d(R,L)$ are isomorphic, for any rank one projective $R$-module $L$.

%	The next lemma asserts that it is enough to prove
%the projective generation of a locally complete intersection ideal (of top height) in $R$ for reduced rings only.  \cite[Proposition 2.15]{MKDSMB}.
%
%	\bl\label{nil}
%	Let $I\subset R$ be such that $\hh(I)=\mu(I/I^2)=d$ and let $\eta$ be the nilradical of $R$. Moreover, assume that there exists a projective $R/\eta$-module $P'$ (with trivial determinant) and a surjection  $\widetilde{ \alpha}:P'\surj \ol{I}$, where `bar' denotes going modulo $\eta$. Then there exists a projective $R$-module $P$ (with trivial determinant) and a surjection $\alpha: P\surj I$. 
%	\el
%	\proof 	Let $R_{\text{red}}=R/\eta$. Let us choose an isomorphism $\chi':R_{\text{red}}\iso \wedge^dP'$ such that $e(P',\chi')=(\ol I)$ in $\E^d(R_{\text{red}})$. As the canonically map $\phi:\E^d(R)\iso \E^d(R_{\text{red}})$ is isomorphic we have $\phi(I)=(\ol I)$ in $\E^d(R_{\text{red}})$. Now applying \cite[Proposition 2.15]{MKDSMB} we get a projective $R$-module $P$ (with trivial determinant) of rank $d$ and an
%	isomorphism $\chi: R \cong  \wedge ^d P$ such that $e(P, \chi) = (I)$ in $\E^d(R)$. Now applying \cite[Corollary 4.3]{SMBB3} we get the required surjection.\qed

The next proposition is due to M. P. Murthy \cite[Theorem 3.3]{M}. For smooth real affine varieties this was proved in \cite[Lemma 5.1]{BRS}.
%\bd
%Let $A$ be a ring of dimension $d$, and let $I\subset A$ be an ideal. An ideal $J \in A$ is said to be residual to $I$ if
%\begin{enumerate}[\quad\quad\quad (1)]
%	\item  $IJ$ is generated by $d$ elements,
%	\item  $I+J=A$,
%	\item $J$ is a local complete intersection of height $d$. 
%\end{enumerate}
%\ed

	\bp\label{PG2}
	In addition to condition \textbf{P}, let $R$ be reduced. Let $J\subset R$ be a reduced ideal such that $\mu(J/J^2)=\hh(J)=d$. Let $L$ be a rank one projective $R$-module. Then there exist (a)  an ideal $K\subset R$ of height $d$, which is a product of distinct smooth complex maximal ideals, and (b) a projective $R$-module $P$ of rank $d$ (which is unique up to isomorphism if we fix $K$) with determinant $L$ such that the following hold.
	\begin{enumerate}[\quad\quad\quad (1)]
		\item $J+K=R$,
		\item $J$ is a surjective image of $P$,
		\item $[P]-[ L\oplus R^{d-1}]=[R/K]\in \mathrm{K}_0(R)$,
		\item  $e(P)=(J)$ in $\E^d_0(R,L)$,
	\item 	$(J)+(d-1)!(K)=0 \text{ in } \E^d_0(R,L)$, and additionally, if $J$ is a product of distinct smooth maximal ideals in $R$ of height $d$, then $[J]+(d-1)![K]=0$ in $\text{CH}_0(\Spec(R))$.
	
	\end{enumerate}

	\ep
	\proof
	
	Applying Lemma \ref{SB} we can find an ideal $I\subset R$, comaximal with $J$ such that $I$ is a product of finitely many distinct smooth complex maximal ideals of height $d$, and $(IJ)=(I)+(J)=0$ in $\E_0^d(R,L)$.
	
	Since $\E^d_0(R,L)$ is a divisible group by Proposition \ref{decg}, we get an ideal $K\subset R$ with $\hh(K)=\mu(K/K^2)=d$ such that $(I)=(d-1)!(K)$ in $\E^d_0(R,L)$. Moreover, using Lemma \ref{SB} twice we may choose $K$ such that $K+JI=R$, and $K$ is a product of distinct smooth complex maximal ideals of height $d$. Let $K=\langle a_1,\ldots,a_d\rangle +K^2$, and let $K'=\langle a_1,\ldots,a_{d-1}\rangle +K^{(d-1)!}$ be the Boratynski ideal of $K$. Since $K$ is comaximal with $JI$, it follows from the definition of $K'$ that $K'+JI=R$. Then by Lemma \ref{BDML} we obtain that $(K')=(d-1)!(K)$ in $\E^d_0(R,L)$. Now as $(I)+(J)=0$, we get the following.
	\begin{equation}\tag{$\star$}
	(J)+(d-1)!(K)=0 \text{ in } \E^d_0(R,L)
	\end{equation} Hence $(JK')=(J)+(K')=0$ in $\E^d_0(R,L)$. Therefore, we get $(I)=(K')$ in $\E_0^d(R,L)$. By applying \cite[Theorem 2.2]{M} we get a (unique up to isomorphism if we fix $K$) projective $R$-module $Q$ of rank $d$ with trivial determinant such that (i) $K'$ is a surjective image of $Q$ and (ii) $[Q]-[R^{d}]=-[R/K]$ in $\mathrm{K}_0(R)$. We point out that the uniqueness of $Q$ such that $Q$ satisfies (i) and (ii) follows from Theorem \ref{PC}. Note that, since $(JK')=0$ in $\E^d_0(R,L)$, there is a surjection $L\oplus R^{d-1}\surj JK'$. Now applying \cite[Theorem 1.3 and Remark 1.4]{M} there exists a (unique up to isomorphism if we fix $K$) projective $R$-module $P$ of rank $d$ and a surjection $P\surj J$ such that $P\oplus Q\cong R^{d}\oplus( L\oplus R^{d-1})$. One may observe that $\det(P)\cong L$. Now let $z=[P]-[L\oplus R^{d-1}]\in \mathrm{K}_0(R)$. Then we obtain the following.
	$$z=[P]-[ L\oplus R^{d-1}]=[R^{d}]-[Q]=[R/K]\in \mathrm{K}_0(R)$$
Since $J$ is a surjective image of $P$, by definition \cite[$\S 6$]{SMBB3} we get (4).

Now we assume that $J$ is a product of distinct smooth maximal ideals in $R$ of height $d$. Since there exists a canonical homomorphism $$\E^d_0(R,L)\to \text{CH}_0(\Spec(R))\text{ sending } (m)\mapsto [m],$$ for every smooth maximal ideals in $R$, (5) follows from $(\star)$. This concludes the proof.\qed

	Now we are ready to prove the main theorem of this section.
	 
\bt\label{QBRS}
In addition to condition \textbf{P}, let $R$ be reduced. Let $J\subset R$ be an ideal which is a product of distinct smooth maximal ideals in $R$ of height $d$. If the cycle $[J]$ associated to $R/J$ is $0$ in $\text{CH}_0(\Spec(R))$, then $\mu(J)=d$. Furthermore, the canonical map $\theta_R:\E^d_0(R)\iso \text{CH}_0(\Spec(R))$ is an isomorphism.
\et
	 
\proof Let $X=\Spec(R)$. Applying Proposition \ref{PG2}, there exists a projective $R$-module $P$ of rank $d$ with trivial determinant and a reduced ideal $K$ of height $d$, which is a product of distinct smooth complex maximal ideals such that the following hold.
\begin{enumerate}[\quad\quad\quad (1)]
	\item $J+K=R$,
	\item $J$ is a surjective image of $P$,
	\item $[P]-[R^{d}]=[R/K]\in \mathrm{K}_0(R)$,
	\item  $e(P)=(J)$ in $\E^d_0(R)$,
	\item $[J]+(d-1)![K]=0$ in $\text{CH}_0(X)$.
\end{enumerate}
Since $[J]=0$ in $\text{CH}_0(X)$, from (5) we get the following. \begin{equation}\tag{$\star$}
	(d-1)![K]=0 \text{ in }  \text{CH}_0(X)
\end{equation}
 Let us define $R_{\mathbb C} := \frac{R[T]}{\langle T^2+1\rangle }\cong R\otimes_{\mathbb R}\BC$, and let $X_{\BC}=\Spec(R_{\BC})$. As $K$ is supported by only smooth complex maximal ideals, the cycle $[K]\in\text{CH}_0(X)$ sits inside the image of the canonical homomorphism $\text{CH}_0(X_{\BC})\to \text{CH}_0(X)$. Then, since $\text{CH}_0(X_{\BC})$ is divisible by \cite[Theorem 6.7]{AKMU}, it follows that $[K]\in \text{D}(X)$, where $\text{D}(X)$ is the maximal divisible subgroup of $\text{CH}_0(X)$. Moreover, from ($\star$) we get that $[K]$ is a torsion element in $ \text{D}(X)$. It follows from \cite[Theorem 1.3]{AKMU} that $\text{CH}_0(X_{\BC})$ is torsion free. Therefore, applying \cite[Proposition 1.4 (b)]{JLCTCS} we obtain that the torsion subgroup  $ \text{D}(X)_{\text{tor}}$ is $0$. In particular, this will give us $[K]=0$ in $\text{CH}_0(X)$. We recall that there is a natural homomorphism $\text{CH}_0(X)\to \mathrm{K}_0(R)$ sending $[K]\mapsto [R/K]$ \cite[$\S 2$]{LW}. Therefore, from (3) we get that $P$ is a stably free module of rank $d$. Hence, using Theorem \ref{P} one may obtain that $P$ is free, in particular $\mu(J)=d$. 

Now we shall show that the canonical map $\theta_R:\E^d_0(R)\iso \text{CH}_0(X)$ is an isomorphism. Note that we have already proved the injective part. To elaborate this let $I\subset R$ be an ideal such that $\mu(I/I^2)=\hh(I)=d$ and $\theta_R(I)=0$ in $ \text{CH}_0(X)$. Applying Lemma \ref{SB} twice we may assume that $I$ is a product of distinct smooth maximal ideals in $R$ of height $d$. Then $\mu(I)=d$ follows from the previous paragraph. 

Since any smooth maximal ideal $m\subset R$ of height $d$ satisfies the condition $\mu(m/m^2)=d=\hh(m)$, the surjectivity of the map $\theta_R$ follows from the fact that $(m)\in \E^d_0(R)$. This concludes the proof. \qed

%
%\bd
%Let $A$ be a reduced affine algebra of dimension $d$ over an infinite perfect field. Let $P$ be a projective $A$-module of rank $d$. The d-th Chern class of $P$ is defined by $$C_d(P):= \sum_{i=0}^d \wedge^i (P*).$$
%\ed
%
%\rmk It follows from \cite[Remark 3.6]{M} that there exists a locally complete intersection reduced ideal $I\subset A$ of height $d$ such that (1) $I$ is a product of distinct smooth maximal ideals, and (2) $C_d(P)=[A/I]\in \text{F}^d\mathrm{K}_0(A)$.
% 

\subsection{Consequences}\label{1C} Now we obtain the following series of results as a consequence of the results proved in this section until now. We commence with an analogy to A. A. Rojtman's theorem \cite{Roj}.

\bc\label{tf}
The $d$-th Euler class group $\E^d(R)$ is torsion free.
\ec

\proof By \cite[Corollary 4.6]{SMBB3} it is enough to assume that $R$ is reduced. Let $I\subset R$ be an ideal with $\mu(I/I^2)=\hh(I)=d$ such that $n(I)=0$ in $\E_0^d(R)$, for some $n\ge 1$. Using Lemma \ref{SB} twice we may assume that $I$ is a product of finitely many distinct smooth complex maximal ideals of height $d$. It follows from Corollary \ref{Cor} and Theorem \ref{QBRS} that it is enough to show $[I]=0$ in $\text{CH}_0(\Spec(R))$. Since $I$ is a product of smooth complex maximal ideals, $[I]$ sits inside the image of the canonical map $\text{CH}_0(\Spec(R_{\BC}))\to \text{CH}_0(\Spec(R))$. Hence, by \cite[Proposition 1.4 (b)]{JLCTCS} and \cite[ Theorem 1.3]{AKMU} it follows that $[I]=0$ in $\text{CH}_0(\Spec(R))$. \qed

\smallskip

%\rmk One may deduce Corollary \ref{tf} from Remark \ref{inv} as well.

\rmk\label{unique} Since $\E^d_0(R,L)$ is torsion free, in Proposition \ref{PG2} the choice of the pair $(K,P)$ is unique in the sense that if there exists another pair $(K'',P'')$ satisfying all the conditions of \ref{PG2}, then $(K)=(K'')$ in $\E^d_0(R)$ and $P\cong P''$.

 Since  $R\inj R_{\BC}$ is an integral extension, it induces natural homomorphisms $\chi_R:\E^d(R)\to \E^d(R_{\BC})$ and $\gamma_R:\E_0^d(R)\to \E_0^d(R_{\BC})$ (see \cite[Definition 3.2]{YY}).

\bc\label{pfpb}
The natural map $\chi_R:\E^d(R)\to \E^d(R_{\BC})$ is injective.
\ec

\proof Without loss of generality we may assume that $R$ is reduced. Let us take $X=\Spec(R)$, and $X_{\BC}=\Spec(R_{\BC})$. We note that, the base change map $\pi:X_{\BC}\to X$ is a finite locally free morphism of degree $2$. Hence, the map $\pi$ will canonically induce homomorphisms $\pi^*:\text{CH}_0(X)\to \text{CH}_0(X_{\BC})$ and $\pi_*:\text{CH}_0(X_{\BC})\to \text{CH}_0(X)$ such that the composition $\pi_*\circ \pi^*:\text{CH}_0(X)\to \text{CH}_0(X)$ is multiplication by $2$. Since $\text{CH}_0(X)$ is uniquely divisible by Corollary \ref{tf}, the composition $\pi_* \circ \pi^*$ is an isomorphism. Therefore, the map $\pi^*$ is injective. Now, we consider the following diagram, in which each square is commutative.
$$\begin{tikzcd}
		\E^d(R) \arrow{r}{\chi_R} \arrow{d}{\Phi_R} & \E^d(R_{\BC})\arrow{d}{\Phi_{R_{\BC}}} \\
	\E^d_0(R) \arrow{r}{\gamma_R } \arrow{d}{\theta_R} & \E^d_0(R_{\BC})\arrow{d}{\theta_{R_{\BC}}}\\
	\text{CH}_0(X)\arrow{r}{\pi^*} & \text{CH}_0(X_{\BC})
\end{tikzcd}$$One may observe that it follows from Corollary \ref{Cor}, Theorem \ref{QBRS}, and \cite{AKMU} that all the vertical maps are isomorphisms. Hence $\chi_R$ is injective. This concludes the proof.
  \qed

%We claim that $\Gamma_R$ is surjective. To see the surjectivity of $\Gamma_R$ we conclude as follows: by Step - 3, it is enough to find preimages of all smooth complex maximal ideals in $R$ of height $d$. Let $m\subset R$ be a smooth complex maximal ideal of height $d$, then we have $\mu(m/m^2)=d=\hh(m)$. It follows from the previous discussion that the extension of $m$ in $R_{\BC}$ (say $\ol m$) is also a maximal ideal in $R_{\BC}$ of height $d$. Since $\mu(m/m^2)=d$, we have $\mu(\ol m/\ol m^2)=d$. Therefore, we get $(\ol m)\in \E^d_0(R_{\BC})$. Now it was proved in \cite[Lemma 7.4]{AKMU} that the groups $\E^d_0(R_{\BC})$ and $\E_s(R_{\BC})$ are canonically isomorphic. This proves that $\Gamma_R$ is surjective. 

\bc
Let $P$ be as in Proposition \ref{PG2}. Then the weak Euler class $e(P)=0$ in $\E^d_0(R,L)$ if and only if $P\cong L\oplus R^{d-1}$. 
\ec
\proof First, we assume that $e(P)=0$ in $\E^d_0(R,L)$. Then by Proposition \ref{PG2} (4) we get that $J$ is a complete intersection ideal of height $d$. Applying Proposition \ref{PG2} (5) and Corollary \ref{tf} we obtain that $[K]=0$ in $\text{CH}_0(\Spec(R))$. Hence, using Theorem \ref{QBRS} it follows that $(K)=0$ in $\E^d_0(R,L)$. Then using Corollary \ref{Cor} and \cite[Theorem 4.2]{SMBB3} it follows that $K$ is a complete intersection ideal of height $d$, in particular $[R/K]=0$ in $\text{K}_0(R)$. Therefore, from Proposition \ref{PG2} (3) we get that $P$ is a stably isomorphic with $ L\oplus R^{d-1}$. Now it follows from Theorem \ref{PC} that $P\cong  L\oplus R^{d-1}$. The converse part is trivial.  \qed

%	 In summary we get the following: 
%	
%	\begin{enumerate}[\quad\quad\quad (a)]
%		\item $K+J=R$
%		\item a surjection $\beta: R^d\surj JK'$
%	\end{enumerate}
%	
%	Then applying \cite[Theorem 3.1]{} there exists a projective $R$-module $P$ of rank $d$ such that (1), (2) and (3) are satisfies. 
%	 \qed

%	Since there are canonical isomorphism between $F^dK_0(R)\cong E_0(R)\cong E(R)$, we must have $(R/I)=(R/I')$ in $F^dK_0(R)$. Since $J'$ is a surjective image of $P$, in $K_0(R)$ we have $C_d(P^*)=(R/I'),$ and hence $C_d(P^*)=(R/I)$ where $C_d(-)$ denotes the $d$-th Chern class.
%	\par Consider a surjection $\o\alpha:P/IP\surj I/I^2$. Let $f:P\to I$ be any lift of $\ol \alpha $ (i.e. $ f\otimes R/I=\ol \alpha$), then we have $I=f(P)+I^2$. Then we can find $ I''\subset R$ comaximal with $I$ of height $d$ such that $I\cap I''=f(P)$. Thus in $K_0(R)$ we have $(R/I')=(R/I)=C_d(P*)=(R/II'')=(R/I)+(R/I'')$. Thus we get $(R/I'')=0$ in $F^dK_0(R)$ and hence $(I'')=0$ in $E(R)$.
%	\par Let `bar' denotes going modulo $I''$. Get $\delta^{-1} :P/I''P\iso (R/I'')^d$ such that $\wedge^d \delta =\ol \chi$, where $\chi:\wedge^dP\cong R$. Let $\ol \beta=\ol f\delta :(R/I'')^d\surj I''/I''^2$. Since $(I'')=0$ get $\beta :R^d\surj I''$ such that $\beta \otimes R/I''=\ol \beta$ and $f:P\surj I\cap I''$. Then note that $(\beta \otimes A/I'')\delta^{-1} = \ol \beta\delta^{-1}=\ol f$. Hence by Subtraction principal (\cite{SMBB3}, Theorem 3.3) there exists $\theta:P\surj I$ such that $\theta \otimes R/I=f\otimes R/I=\ol \alpha$.\qed \\

	One may compare the following result with \cite[Theorem 2.1]{MM}.
	
	\bc\label{Mur2}
	Let $I\subset R$ be an ideal such that $\mu(I/I^2)=\hh(I)=d$. Let $P$ be a projective $R$-module of rank $d$ with determinant $L$. Suppose that $\ol f:P/IP\surj I/I^2$ be a surjective map. Then there exists a surjective lift  $f:P\surj I$ of $\ol f$ if and only if $e(P)=(I)$ in $\E^d_0(R,L).$ 
	\ec
	\proof First we prove the theorem with the assumption that $R$ is reduced. We fix an isomorphism $\chi:\wedge^dP\cong L$. Note that if $f:P\surj I$ then from the definition (see \cite[Section 6, paragraph 7]{SMBB3}) it follows that $e(P)=(I)$ in $\E^d_0(R,L).$ Therefore, we assume that $e(P)=(I)$. Let $\alpha:P\to I$ be any lift (might not be surjective) of $\ol f $ (i.e. $ \alpha\otimes R/I=\ol f$). Then we get $I=\alpha(P)+I^2$. By \cite[Lemma 2.11]{SMBB3} there exists an ideal $ J\subset R$ co-maximal with $I$ of height $d$ such that $I\cap J=\alpha(P)$. Therefore, in $\E^d_0(R,L)$ we have $(I)=e(P)=(I\cap J)=(I)+(J)$. That is $(J)=0$ in $\E^d_0(R,L)$. Since $\E^d(R,L)\iso \E^d_0(R,L)$, we obtain that any local $L$-orientation $\omega_J$ of $J$ is a global $L$-orientation.

	Let `bar' denote going modulo $J$. We can find an isomorphism $\delta^{-1} :P/JP\iso (L/JL)\oplus (R/J)^{d-1}$ such that $\wedge^d \delta =\ol \chi$. Let us define $\ol \beta:=\ol f\delta :(L/JL)\oplus (R/J)^{d-1}\surj J/J^2$. Since $(J)=0$ in $\E^d_0(R,L)$ by \cite[Corollary 4.4]{SMBB3} and Corollary \ref{Cor} we get a surjection $\beta :L\oplus R^{d-1}\surj J$ such that $\beta \otimes R/J=\ol \beta$. Then we note that $(\beta \otimes R/J)\delta^{-1} = \ol \beta\delta^{-1}=\ol f$. Hence, by Subtraction principal \cite[Theorem 3.3]{SMBB3}  there exists $f:P\surj I$ such that $f \otimes R/I=\alpha\otimes R/I=\ol f$. This concludes the proof when $R$ is reduced.
	
	Now let $\eta$ be the nilradical of $R$. Let $A=R/\eta$, and `tilde' denote going modulo $\eta$. Let us assume that $e(P)=(I)$ in $\E_0^d(R,L)$. From the previous step, there exists a surjective lift $f':\widetilde P\surj \wt I$ of $\ol f$. We choose a lift (might not be surjective) $g:P\to I$ of $f'$. Then we get that $g(P)+I\cap \eta=I$. Now one can obtain the required surjection $f:P\surj I$ using the arguments given in \cite[Corollary 4.6]{SMBB3} or see \cite[Corollary 4.13]{MPHIL} for a detailed version. \qed
%	
%\bl
%There exists a canonical group homomorphism $\alpha:\E^d(R,L)\to \E^d(R)$.
%\el
%
%\proof Hint \cite[6.8]{SMBB3}
	%\proof The if part is follows from the definition of weak Euler class group so without loss of generality we may assume $e(P)=(I)$ in $E_0(R)$. Thus we only need to show that $f$ is a surjective lift of $\ol f$. Following the proof of Corollary \ref{Mur2} get  $ I''\subset R$ co-maximal with $I$ of height $d$ such that $I\cap I''=\alpha(P)$. Thus we have $e(P)=(I\cap I'')=(I)+(I'')=e(P)+(I'')$ that is $(I'')=0$
	
	The next result is an analogy of \cite[Theorem 1]{NMK84} in our set-up.
	
	\bc\label{NMKL}
		Let $I\subset R$ be an ideal such that $\mu(I/I^2)=\hh(I)=d$. Let $P$ be a projective $R$-module of rank $d$ such that $f:P\surj I $ is a surjection. Then $P$ has a unimodular element if and only if $\mu(I)=d$.
	\ec
	\proof Let $L$ be the determinant of $P$, and let $\chi:\wedge^d P\cong L$ be an isomorphism. By Theorem \ref{Mur2} we obtain that $e(P)=(I)$ in $\E^d_0(R,L)$. Let $\{\chi,f\}$ induce the local $L$-orientation $\omega_I$ of $I$. Since the groups $\E^d_0(R,L)$ and $\E_0^d(R)$ are canonically isomorphic \cite[Theorem 6.8]{SMBB3} we obtain that $(I)=0$ in $\E^d_0(R)$ if and only if $(I)=0$ in $\E^d_0(R,L)$. Applying Corollary \ref{Cor} we obtain that $(I)=0$ in $\E^d_0(R,L)$ if and only if $(I,\omega_I)=0$ in $\E^d(R,L)$. Now using \cite[Theorem 4.2]{SMBB3} we get $\mu(I)=d$ if and only if $(I,\omega_I)=0$ in $\E^d(R,L)$. By \cite[Corollary 4.4]{SMBB3} we obtain that $(I,\omega_I)=0$ in $\E^d(R,L)$ if and only if $P$ has a unimodular element. This concludes the proof. \qed

	 %Since the groups $\E^d_0(R,L)$ and $\E_0^d(R)$ are canonically isomorphic \cite[Theorem 6.8]{SMBB3} we obtain that $(I)=0$ in $\E^d_0(R)$.
%	 First we assume that $P$ has a unimodular element. Then it follows from \cite[Corollary 4.4]{SMBB3} that $(I)=0$ in $\E^d_0(R,L)$.  Therefore, it follows from [Proposition \ref{decg}, Step - 1] that $\mu(I)=d$. 

%	\smallskip
%	
%\rmk Let $A$ be as in Corollary \ref{PGM} (1) or (3). 
%	
%	\smallskip
	
%	Let $I=\langle a_1,\cdots,a_d\rangle +I^2$. By (\cite{SMBB3}, Lemma 2.11) get $e\in I^2$ be such that $e(1-e)\in \langle a_1,\cdots,a_d\rangle $. Moreover by a theorem of Eisenbud and Evans (see \cite{SMBB3}, Corollary 2.13) replacing $a_i$ by $a_i+\lambda_ie$ we may assume $\hh(\langle a_1,\cdots,a_d\rangle _e)=d$. Let $J=\langle a_1,\cdots,a_d,1-e\rangle $. Then note that\\
%	$(i)\,\hh(J)=d$,\\
%	$(ii)\, I+J=R$,\\
%	$(iii)\, J=\langle a_1,\cdots,a_d\rangle +J^2$.\\
%	Then by Theorem \ref{PG2} $J$ is projective generated and so is $I$ by Lemma \ref{MPGL}. \qed

We now move towards the polynomial extension of $R$ and prove some results as a corollary of our main theorems in this section. We begin with a theorem on projective generation of a locally complete intersection in $R[T]$ of height $d$. For affine algebras over an algebraically closed field of characteristic $0$ a similar result was proved in \cite{BRS98}. 
	
	\bc\label{PG3}
	Let $I\subset R[T]$ be a locally complete intersection ideal such that $\hh(I)=\mu(I/I^2)=d\ge 3$. Let $\omega_I$ be a local orientation of $I$. Then there exists a projective $R[T]$-module $P$ and an isomorphism $\chi:\wedge^d P\cong R[T]$ such that $e(P,\chi)=(I,\omega_I)$ in $\E^d(R[T])$. As a consequence, we obtain a surjection $\alpha:P\surj I$ such that $(I, \omega_I)$ is obtained from $(\alpha, \chi)$.
	\ec
	\proof Since $\mathbb Q\subset R$ we may assume that there exists an element $\lambda\in R$ such that either $I(\lambda)=R$ or $\hh(I(\lambda ))=d$. Moreover, taking the transformation $T\mapsto T-\lambda$ we may further assume that either $I(0)=R$ or $\hh(I(0))=d$. Then we note that $(I(0),\omega_I(0))$ is in $\E^d(R)$. Applying Proposition \ref{PG2} (3) and Corollary \ref{Cor} there exists projective $R$-module $Q$ together with an isomorphism $\chi':R\cong \wedge^d Q$ such that $e(Q,\chi')=(I(0),\omega_I(0))$ in $\E^d(R)$. Now using \cite[Theorem 3.5]{MKDSMB} there exists a projective $R[T]$-module $P$ and an isomorphism $\chi:\wedge^d P\cong R[T]$ such that $e(P,\chi)=(I,\omega_I)$ in $\E^d(R[T])$. It follows from \cite[Corollary 4.10]{MKD1} that there exists a surjection $\alpha:P\surj I$ such that $(I, \omega_I)$ is obtained from $(\alpha, \chi)$.\qed

%	Let $\sigma: Q[T]/IQ[T]\cong (R[T]/I)^{d-1}\oplus (L/IL)$ be an isomorphism. We define a surjection $\ol \omega:=\omega_I\sigma: Q[T]/IQ[T]\surj I/I^2$. This will induce $\ol{\omega(0)}:Q/I(0)Q\surj I(0)/I(0)^2$. By Corollary \ref{Mur2} we can lift $\ol{\omega(0)}$ surjectively. Altering $f$ with a surjective lift of $\ol{\omega(0)}$, we may further assume $f\otimes R/I(0)=\ol{\omega(0)}$. 
%	
%	
%	
%	
%	
%	Then by \cite[Remark 3.9]{BR} $\ol \omega$ can be lifted to  $\omega(T):P[T]\surj I/I^2T$. Using \cite[Theorem 3.6]{MKDSMB} there exists a projective $R[T]$ module $Q$ of rank $d$ with $Q/TQ\cong P$ such that $Q\surj I$. This concludes the proof when $\hh(I(0))=d$. 
%	
%	Now, if $I(0)=R$ then by \cite[Remark 3.9]{BR} any local orientation $\omega$ of $I$ can be lifted to a surjection $(R[T])^d\surj I/I^2T$. Then again as before applying \cite[Theorem 3.6]{MKDSMB} the proof follows.\qed 
	
	The following corollary gives an affirmative answer to a question asked in \cite[Question 2]{MKD1} in our set-up.
	
	\bc
	Let $R(T)$ be the ring obtained from $R[T]$ by inverting all monic polynomials in $R[T]$. Then for all $d\ge 3$ the canonical map $$\Gamma :\E^d(R[T])\to \E^d(R(T))$$is injective.
	\ec
	\proof Let $(I,\omega_I)$ in $\E^d(R[T])$ be such that $(IR(T),\omega_I\otimes R(T))=0$ in $\E^d(R(T))$. Applying Proposition \ref{PG3} there exists a projective $R[T]$-module $P$ together with an isomorphism $\chi:\wedge^d P\cong R[T]$, and a surjection $\alpha:P\surj I$ such that (1) $e(P,\chi)=(I,\omega_I)$ in $\E^d(R[T])$, and (2) $(I,\omega_I)$ is obtained from $(\alpha,\chi)$. Since $(IR(T),\omega_I\otimes R(T))=0$ by \cite[Corollary 4.4]{SMBB3} we obtain that $P\otimes R(T)$ has a unimodular element. Hence, applying \cite[Theorem 3.4]{BRS01} we get that $P$ has a unimodular element. Therefore, using \cite[Corollary 4.11]{MKD1} it follows that $(I,\omega_I)=0$ in $\E^d(R[T])$. This concludes the proof.\qed
	
\section{A nice group structure of $\text{WMS}_{d+1}(R)$}\label{2}
In this section, we prove that $\text{WMS}_{d+1}(R)$ has a nice group structure. Hence, it follows from \cite{VdKM} that $\text{WMS}_{d+1}(R)$ is canonically isomorphic to the universal Mennicke symbol  $\text{MS}_{d+1}(R)$. The outline of the proof follows the arguments given in \cite[Theorem 3.9]{AR}.

\bt\label{NGS}
The abelian group $\text{WMS}_{d+1}(R)$ has a nice group structure. That is for any $(a,a_1,\ldots,a_d)$ and $(b,a_1,\ldots,a_d)\in \Um_{d+1}(R)$ we have the following identity: $$[(a,a_1,\ldots,a_d)]\star[(b,a_1,\ldots,a_d)]=[(ab,a_1,\ldots,a_d)].$$ In particular $\text{WMS}_{d+1}(R)\cong \text{MS}_{d+1}(R)$.
\et
\proof Without loss of generality we may assume that $R$ is a reduced ring (see \cite[Lemma 3.5]{AR}). Moreover, if $R$ satisfies condition \textbf{P}-2, then by Lemma \ref{PR} taking $\BS$ to be the collection of all real maximal ideals we may further assume that for any $2\le i\le d$, the ring $R/\langle a_i,a_{i+1},\ldots,a_d\rangle $ is a smooth real affine algebra of dimension $i-1$, having no real maximal deal. On the other hand, if $R$ satisfies \textbf{P}-1, then using Theorem \ref{SwB} one can achieve the same. Applying the product formula as given in \cite{VdK} we get $$[( a, a_1,\ldots,a_d)]\star[( b, a_1,\ldots,a_d)]=[(a(b+p)-1,(b+p)a_1,a_2,\ldots,a_d)],$$
where $p$ is chosen such that $ap-1\in \langle a_2,\ldots,a_d\rangle $. Let $B=R/\langle a_2,a_{3},\ldots,a_d\rangle $ and let `bar' denote going modulo $\langle a_2,a_{3},\ldots,a_d\rangle $. Then by the Bass-Kubota theorem \cite[Theorem 2.12]{VdK} we have $\mathrm{SK_1}(B)\cong \text{MS}_2(B)$. Therefore, in the group $\text{MS}_2(B)$ we get the following.
$$[(\ol a(\ol b+\ol p)-\ol 1,(\ol b+\ol p)\ol a_1)]=[(\ol a(\ol b+\ol p)-\ol 1,\ol a_1)]$$
Hence there exists a $\sigma\in \SL_2(B)\cap \E_3(B)$ such that: $$(\ol a(\ol b+\ol p)-\ol 1,(\ol b+\ol p)\ol a_1)\sigma =(\ol a(\ol b+\ol p)-\ol 1,\ol a_1).$$ By Theorem \ref{pa} we obtain that $\sigma\in \SL_2(B)\cap \mathrm{ESp}(B)$. Then applying \cite[Corollary 2.3]{AAS} we can find $\alpha\in \E_{d+1}(R)$ such that the following identity holds:$$ (a(b+p)-1,(b+p)a_1,a_2,\ldots,a_d)\alpha=(a(b+p)-1,a_1,a_2,\ldots,a_d).$$
Now as $ap-1\in \langle a_2,\ldots,a_d\rangle $ we get the following identity: $$[(a(b+p)-1,a_1,a_2,\ldots,a_d)]=[(ab,a_1,a_2,\ldots,a_d)].$$
This concludes the proof. \qed

The remaining part of the section is devoted to establishing some consequences of the above theorem.

\bc\label{ngsrc}
Let $I\subset R$ be an ideal. Then the abelian group $\frac{\Um_{d+1}(R,I)}{\E_{d+1}(R,I)}$ has a nice group structure. That is $$[(a,a_1,\ldots,a_d)]\star[(b,a_1,\ldots,a_d)]=[(ab,a_1,\ldots,a_d)],$$ where $[-]$ denote the class in the relative elementary orbit space of relative unimodular rows of length $d+1$.
\ec
\proof Applying Theorem \ref{NGS} and Lemma \ref{EL} it follows that the group $\text{WMS}_{d+1}(R\oplus I)$ has a nice group structure. Hence, by \cite[Lemma 3.6]{AAR} the proof concludes. \qed

We end this section with an extension of a result due to J. Fasel \cite[Theorem 2.2]{Fasel}.

\bc\label{NGSC} The group $\text{WMS}_{d+1}(R)$ is divisible.
\ec
% \ \textbf{Divisibility} 
\proof  Let $v=(v_0,\ldots,v_d)\in \Um_{d+1}(R)$ and let $n\in \BN$. To prove the result we shall find a $w\in \Um_{d+1}(R)$ such that $[v]=[w]^n$. By Proposition \ref{PR} (taking $\BS$ to be the collection of all real maximal ideals, whenever $R$ satisfies \textbf{P}-2) or Theorem \ref{SwB} we may assume that $R/\langle v_2,\ldots,v_d\rangle $ is a smooth curve having no real maximal ideal. Let $C=R/\langle v_2,\ldots,v_d\rangle $, and let `bar' denote going modulo $\langle v_2,\ldots,v_d\rangle $. Then by the Bass-Kubota theorem \cite[Theorem 2.12]{VdK}, we have $\text{MS}_2(C)=\mathrm{SK_1}(C)$. The latter one is a divisible group by Proposition \ref{pa}. Now following the argument given in [Theorem \ref{P}, paragraph 2], one can get a $\gamma\in \E_{d+1}(R)$ such that $v=\gamma(u^n_0,u_1,v_2,\ldots,v_d)$. We define $w:=(u_0,u_1,v_2,\ldots,v_d)$. Then applying Theorem \ref{NGS} we obtain the following identity: 
$$[v]=[(u^n_0,u_1,v_2,\ldots,v_d)]=[(u_0,u_1,v_2,\ldots,v_d)]^n=[w]^n.$$This concludes the proof.\qed

\section{Improved stability for $\mathrm{K}_1$ and $\mathrm{K}_1\mathrm{Sp}$ groups}\label{3}
In addition to condition \textbf{P}, let $R$ be a regular domain and let $I\subset R$ be a principal ideal. In this set-up, we prove that the injective stability of the groups $\SK(R,I)$ and $\Sp(R,I)$ decreases by one. This improves L. N. Vaser{\v{s}}te{\u{\i}}n's general stability bounds \cite{V69}, \cite{Va1}, and \cite{Vas}. The arguments used in this section are motivated from \cite{RV}.

\bt\label{ISK1}
In addition to condition \textbf{P}, let $R$ be an integral domain, and let $I=\langle a\rangle \subset R$ be a principal ideal. Let $\sigma\in \SL_{d+1}(R,I)\cap \E(R,I)$. Then $\sigma$ is isotopic to identity. Moreover, if we assume that $R$ is regular, then $\mathrm{SK_1}(R,I)\cong\frac{\SL_{d+1}(R,I)}{\E_{d+1}(R,I)}$. 
\et
\proof The surjectivity of the canonical map $\Gamma_{d+1}:\frac{\SL_{d+1}(R,I)}{\E_{d+1}(R,I)}\to \mathrm{SK_1}(R,I)$ follows from \cite{V69}. To prove the injectivity, first, we note that the only non-trivial part is to prove the following. $$\E(R,I) \cap \SL_{d+1}(R,I)\subset \E_{d+1}(R,I)$$Let $\sigma\in \E(R,I) \cap \SL_{d+1}(R,I)$. By stability theorem due to  L. N. Vaser{\v{s}}te{\u{\i}}n \cite{V69} we get $\sigma\in \E_{d+2}(R,I)\cap \SL_{d+1}(R,I)$. It follows from \cite[Chapter 1, $\S$ 2]{SV} that $1\perp \sigma$ is of the form $$\prod\E_{ij}(b)\E_{ji}(a\lambda_k)\E_{ij}(-b),$$ where the product runs for some $i\not= j$ and $b,\lambda_k\in R$. We define the following.$$\tau(T):=\prod\E_{ij}(b)\E_{ji}(\lambda_k T)\E_{ij}(-b)$$ Then $\tau(T)\in \E_{d+2}(R[T],\langle T\rangle )$ such that $\tau(a)=1\perp\sigma$. Let $t=T^2-Ta\in R[T]$ be a non-zero divisor, and let $v=e_1\tau(T)$. Then $v\in \Um_{d+2}(R[T],\langle t\rangle )$. Since $R$ is an affine domain, any maximal ideal in $R[T]$ is of height $d+1$. Hence, with the help of Suslin's monic polynomial theorem (see \cite[Chapter III, $\S 3$, 3.3, page no 108]{Lam}) one may observe that $R/\mm\cap R \inj R[T]/\mm$ is an integral extension, for any maximal ideal $\mm\subset R[T]$. This in particular gives us that $R[T]$ satisfies \textbf{P} [cf. the proof of Lemma \ref{EL}]. Therefore, by Proposition \ref{S1} there exists $\chi(T)\in \SL_{d+2}(R[T],\langle t\rangle )$, such that  $v=e_1\chi(T)$. Since $e_1\tau(T)\chi(T)^{-1}=e_1 $, the matrix $\tau(T)\chi(T)^{-1}$ is of the form $$(1\perp \rho(T))\prod_{i=1}^{d+2}\E_{i,1}(\lambda_i),$$where $\lambda_i\in \langle t\rangle $, $\rho(T)\in \SL_{d+1}(R[T],\langle T\rangle )$ and $\rho(a)=\sigma$. Now since $\chi(T)\equiv \text{I}_{d+2}$ modulo $\langle t\rangle $, we have $\chi(0)=\chi(a)=\text{I}_{d+2}$. Hence, $\rho(0)=\text{I}_{d+1}$. In other words we get $\rho(T)\in \SL_{d+1}(R[T],\langle T\rangle )$ is an isotopy of $\sigma.$

Now, if we assume that $R$ is regular, then by \cite[Theorem 3.3]{TV} we get $\rho(T)\in \E_{d+1}(R[T],\langle T\rangle )$. Hence $\sigma=\rho(a)\in \E_{d+1}(R,\langle a \rangle)$. This concludes the proof.\qed

%The following corollary might be seen as an extension of Proposition \ref{po},
%
%\bc
%Additionally, let $R$ be a regular domain. Then $\mathrm{SK_1}(R)$ is a divisible group. 
%\ec
%\proof It was proved in \cite[Theorem 3.16(iv)]{VdK} the first row map $\SL_{d+1}(R)\to \text{WMS}_{d+1}(R)$ sending $\alpha\mapsto [e_1\alpha]$ is a group homomorphism, where $\alpha\in \SL_{d+1}(R)$. Now if $\beta\in \E_{d+1}(R)$, then clearly $[e_1\beta]=0$ in $\text{WMS}_{d+1}(R)$. This 

For the rest of this section, we shift towards the symplectic matrices, and prove an analogous result of Theorem \ref{ISK1} for the symplectic group $\mathrm{K}_1\mathrm{Sp}(R,I)$.

%Using the proof of Lemma \ref{PR} we may replace $v_i$ with $v_i+\lambda_iv_d$ $(i=0,\cdots d-1)$ and make the assumption that $C:=R/\langle v_0,\cdots,v_{d-1}\rangle $ is a finite product of $\mathbb C$. Since $v_d\in I$, this replacements of $v_i$'s does not alter the fact that $(v_0,\cdots,v_d)\in \Um_{d+1}(R,I)$. Let `bar' denote going modulo $\langle v_0,\cdots,v_{d-1}\rangle $. Then there exists $b\in R$ such that the following hold.$$v_d-b^{d!}\in \langle v_0,\cdots,v_{d-1}\rangle$$

\bp\label{RAOVDK}
Let $d\equiv 1\mod(4)$, and let $I=\langle a\rangle \subset R$ be a principal ideal, where $a\in R$ is a non-zero divisor. Then  $$\Um_{d+1}(R,I)=e_1\mathrm{Sp}_{d+1}(R,I).$$
\ep
\proof First, we note that for $d=1$, we have $\mathrm{Sp}_2(R,I)=\SL_2(R,I)$. Hence, there is nothing to prove when $d=1$, so we may assume that $d\ge 5$. Let $v\in \Um_{d+1}(R,I)$.  Recall that, for any unimodular row $w=(w_0,\ldots,w_{r})$, we denote $\chi_{r!}(w)$ by the unimodular row $(w_0,\ldots,w_{r-1}, w_r^{r!})$. It follows from Lemma \ref{bla} that 
\begin{equation*} v\equiv \chi_{d!}(v')\mod \mathrm{E}_{d+1}(R,I),
\end{equation*}
for some $v'\in \Um_{d+1}(R,I)$. Then applying \cite[Theorem 5.5]{PRRA} we obtain the following. \begin{equation} v\equiv \chi_{d!}(v')\mod \mathrm{ESp}_{d+1}(R,I)
\end{equation} By Corollary \ref{ngsrc} it follows that the group $\frac{\Um_{d+1}(R,I)}{\E_{d+1}(R,I)}$ has a nice group structure. Hence, using \cite[Proposition 3.1]{BCR} we get the following. \begin{equation}
 	\chi_{d!}(v')\equiv e_1 \mod \mathrm{Sp}_{d+1}(R,I)
 \end{equation}
Therefore, combing (1) and (2) we get $v\equiv e_1 \mod \mathrm{Sp}_{d+1}(R,I)$. This completes the proof. \qed

\bt\label{ISPS}
In addition to condition \textbf{P}, let $R$ be a regular domain of dimension $d\ge 4$, and let $I=\langle a\rangle \subset R$ be a principal ideal. Moreover, assume that $d \not\equiv 2\mod 4$. Let $n=2[\frac{d+1}{2}]$, where $[\,-\,]$ denotes the smallest integer less than or equals to $-$. Then $\mathrm{Sp}_n(R,I)\cap\mathrm{ESp}(R,I)=\mathrm{ESp}_n(R,I)$. Consequently, the canonical map $\frac{\mathrm{Sp}_n(R,I)}{\mathrm{ESp}_n(R,I)}\to \mathrm{K}_1\mathrm{Sp}(R,I)$ is injective.
\et

\proof First, we observe that applying  L. N. Vaser{\v{s}}te{\u{\i}}n's stability theorem \cite{Vas} it is enough to show that the following. $$\mathrm{Sp}_n(R,I)\cap \mathrm{ESp}_{n+2}(R,I)\subset \mathrm{ESp}_n(R,I)$$ Let $\sigma\in \mathrm{Sp}_{n}(R,I)\cap \mathrm{ESp}_{n+2}(R,I)$. Then $\text{I}_2\perp \sigma$ is of the form $$\prod g\text{s}\E_{ij}(a\lambda_k)g^{-1},$$ where the product runs for some $i\not= j$, $g\in \text{ESp}_{n+2}(R)$ and $\lambda_k\in R$ (cf. \cite[Lemma 3.5, last paragraph]{RPJA}. We define the following. $$\rho(T):=\prod g\text{s}\E_{ji}(\lambda_k T)g^{-1}$$ Then $\rho(T)\in \text{ESp}_{n+2}(R[T],\langle T\rangle )$ such that $\rho(a)=\text{I}_2\perp\sigma$. Let us take $v(T)=e_1\rho(T).$ Then we note that $v(T)\in \Um_{n+2}(R[T],\langle T^2-aT\rangle )$. Since $R$ is an affine domain, any maximal ideal in $R[T]$ is of height $d+1$. Hence, with the help of Suslin's monic polynomial theorem, one may observe that $R/\mm\cap R \inj R[T]/\mm$ is an integral extension, for any maximal ideal $\mm\subset R[T]$. This in particular gives us that $R[T]$ satisfies condition \textbf{P}.

If $d$ is odd, then we have $n=d+1$. Then applying \cite[Lemma 6.3 and Theorem 5.5]{PRRA} we obtain the following. $$v(T)\in \Um_{d+3}(R[T],\langle T^2-aT\rangle )=e_1\E_{d+3}(R[T],\langle T^2-aT\rangle)=e_1\text{ESp}_{d+3}(R[T],\langle T^2-aT\rangle)$$ Hence there exists an $\alpha(T)\in \mathrm{ESp}_{d+3}(R[T],\langle T^2-aT\rangle )$ such that  $v(T)=e_1\alpha(T)$. 

Now if $4$ divides $d$, then we note that $n=d$. Applying Proposition \ref{RAOVDK} one can find an $\alpha(T)\in \mathrm{Sp}_{d+2}(R[T],\langle T^2-aT\rangle )$ such that  $v(T)=e_1\alpha(T)$.

Therefore, in either case there exists an $\alpha(T)\in \mathrm{Sp}_{n+2}(R[T],\langle T^2-aT\rangle )$ such that $v(T)=e_1\alpha(T)$. Then the matrix $\rho(T)\alpha(T)^{-1}$ is of the form $$ \begin{pmatrix}
	1 & 0 & 0\\
	* & 1 & *\\
	* & 0 & \eta(T)
\end{pmatrix},$$
for some $\eta(T)\in \mathrm{Sp}_n(R[T],\langle T\rangle ).$ Since $R$ is regular applying \cite[Theorem 5.3]{AG} we obtain that $\eta(T)\in \mathrm{ESp}_{n}(R[T],\langle T\rangle )$. Since $\alpha(a)=\text{I}_{n+2}$, the matrix $\eta(T)\in \mathrm{ESp}_{n}(R[T],\langle T\rangle )$ must satisfy the equality $\eta(a)=\sigma$. In particular, this gives us $\sigma\in \mathrm{ESp}_{n}(R,I)$. This concludes the proof. \qed

\rmk In addition to condition \textbf{P}, let $R$ be a regular domain. When $d\equiv 3\mod(4)$, one can employ Theorem \ref{ISK1} and follow the argument in \cite[Theorem 3.6]{ST} to conclude that $\Um_{d+1}(R)=e_1\mathrm{Sp}_{d+1}(R)$. However, whether an improvement of the existing symplectic injective stability for $\mathrm{K_1Sp}(R)$ is possible still remains unknown for the case $d\equiv 2\mod (4)$.

\section{A sufficient condition for efficient generation of modules}\label{6}

 This section is motivated from \cite{NMK84} and \cite[$\S 4$]{M}. Here we investigate relations between some Eisenbud-Evans type theorems as studied in \cite{NMK84}. To do so we need some preparation. We begin by recalling the following definition from \cite[Section 1]{Mak}.
 
 \smallskip
 
\bd
Let $A$ be a ring of dimension $d$. The group $\mathrm{F^dK}_0(A)$ is a subgroup of $\mathrm{K_0}(A)$ defined as follows: $$\mathrm{F^dK}_0(A):=\{[A/I]\in \mathrm{K_0(A)}: I\subset A \text{ is a locally complete intersection ideal of height }d\}.$$ 

\ed 
%	This version will be required in Section \ref{}.

Let $A$ be a ring (not necessarily containing $\mathbb Q$). Then one can still define the $d$-th Euler class group $\E^d(A)$ and prove that $(K,\omega_K)=0$ if and only if $\omega_K$ is a global orientation of $K$ (see \cite[Section 4]{BRS5}). Let $J\subset A$ be an ideal such that $\mu(J/J^2)=\hh(J)=d$. Since in a local ring $A_p$, one can always lift generators of $J_p/J_p^2$ (where $p\in V(J)$) to generators of $J_p$, the ideal $J$ is in fact a locally complete intersection ideal of height $d$. Moreover, the group $\mathrm{F^dK}_0(A)$ kills all the complete intersection ideals of height $d$. It is evident that there exists a canonical group homomorphism $\kappa_A:\E^d(A)\to \mathrm{F^dK_0}(A)$ sending $(I,\omega_I)\mapsto [A/I]$, where $\omega_I$ is a local orientation of $I$. In the next theorem, we prove that $\kappa_R$ is an isomorphism.

\bt\label{k_0}
The canonical map $\kappa_R:\E^d(R)\iso \mathrm{F^dK_0}(R)$ is an isomorphism.
\et
\proof Without loss of generality we may assume that $R$ is reduced. Moreover, by Corollary \ref{Cor} it is enough to show that the composition map $\E_0^d(R)\to \mathrm{F^dK_0}(R) $ is an isomorphism. With an abuse of notation, we denote the composition map $\E_0^d(R)\to \mathrm{F^dK_0}(R) $ by the same notation $\kappa_R$. 

In our set-up, the injectivity of $\kappa_R$ is a consequence of results due to A. A. Suslin, M. Boratynski and M. P. Murthy. We elaborate: let $J\subset R$ be an ideal such that $\hh(J)=\mu(J/J^2)=d$, and $(J)\in \ker(\kappa_R)$. Let $J=\langle a_1,\ldots, a_d\rangle +J^2$. We consider the Boratynski ideal $I:=\langle a_1,\ldots, a_{d-1}\rangle +J^{(d-1)!}$ of $J$. Then by \cite[Theorem 2.2]{M} there exists a projective $R$-module $P$ of rank $d$ such that (1) $I$ is a surjective image of $P$, and (2) $(P)-(R^d)=-[R/J]\in \mathrm{F^dK_0}(R)$. Moreover, from Lemma \ref{BDML} we get $(I)=(d-1)!(J)$ in $\E_0^d(R).$ Now since $(J)\in \ker(\kappa_R)$ from (2) one can observe that $P$ is stably free. Hence $P$ is free by Theorem \ref{P} and therefore, we obtain that $(d-1)!(J)=(I)=0$ in $\E_0^d(R)$. Since $\E^d_0(R)$ is torsion free by Corollary \ref{tf} it follows that $(J)=0$ in $\E_0^d(R)$.

Now we prove the surjectivity of the map $\kappa_R$. Applying \cite[Theorem 1.5]{Mak} we obtain that $\mathrm{F^dK_0}(R)$ is generated by the classes $[R/m]$, where $m$ runs
through all the smooth maximal ideals in $R$ of height $d$. Since any smooth maximal ideal $m\subset R$ of height $d$ will satisfy $\mu(m/m^2)=d=\hh(m)$, the ideal $m$ has a preimage in $\E^d_0(R)$. This concludes the proof. \qed

 \smallskip

\rmk Let $R$ be reduced, and let $P$ be a projective $R$-module of rank $d$. Corollary \ref{tf} and Theorem \ref{k_0} immediately establish that the top Chern class $\text{c}_d(P)$ of $P$, as defined in \cite{M}, governs the splitting behavior for $P$ (see \cite[Theorem 3.7]{M}).

The next result is a stronger version of Proposition \ref{PG2}. The proof essentially uses Proposition \ref{PG2} and a subtraction principle due to N. Mohan Kumar and M. P. Murthy \cite[Theorem 1.3]{M}. 

\bp\label{PGM}

Let $A$ be a reduced real affine algebra of dimension $d\ge 2$, and let $J\subset A$ be an ideal which is not contained in any minimal prime ideals of $A$ such that $\mu(J/J^2)=d$. Additionally, assume that $A$ satisfies one of the following conditions.
\begin{enumerate}[\quad\quad (i)]
	\item $A$ satisfies \textbf{P}-1;
%	\item $A$ satisfies P-2 and $\hh(J)\ge 2$;
	\item the intersection of all real maximal ideals in $A$ has height of at least $2$.
\end{enumerate}
 Then there exist (a) a reduced ideal $K\subset A$ of height $d$, which is a product of distinct smooth complex maximal ideals, and (b) a projective $A$-module $P$ of rank $d$ with trivial determinant such that the following hold.
 \begin{enumerate}[\quad\quad\quad (1)]
% 	\item $J+K=R$,
 	\item $J$ is a surjective image of $P$,
 	\item $[P]-[A^{d}]=-[A/K]\in \mathrm{K}_0(A)$.
% 	\item  $e(P)=(J)$ in $\E^d_0(A)$.
 %	\item $[J]+(d-1)![K]=0$ in $\text{CH}_0(\Spec(A))$.
 \end{enumerate}
%Moreover $K$ can be chosen in such a way that $K$ is unique in the sense of Remark \ref{unique}, and hence $P$ is unique (up to isomorphism) satisfying (1) and (2). 

\ep
\proof By Lemma \ref{SB}, one can obtain an ideal $I\subset A$ such that 
\begin{enumerate}[\quad \quad (a)]
	\item $I$ is a product of distinct smooth complex maximal ideals of height $d$,
	\item $I+J=A$, and
	\item $\mu(I\cap J)=d$.
\end{enumerate}
Now applying Proposition \ref{PG2} there exist (i) a reduced ideal $K\subset A$ of height $d$, which is a product of distinct smooth complex maximal ideals and (ii) a projective $A$-module $Q$ of rank $d$ with trivial determinant and such that the following hold:

\begin{enumerate}[\quad\quad\quad (I)]
 	\item $I+K=A$,
	\item $I$ is a surjective image of $Q$,
	\item $[Q]-[A^{d}]=[A/K]\in \mathrm{K}_0(A)$,
	\item  $e(Q)=(I)$ in $\E^d_0(A)$,
	\item $[I]+(d-1)![K]=0$ in $\text{CH}_0(\Spec(A))$.
\end{enumerate}
We now show that the hypotheses of \cite[Theorem 1.3]{M} are satisfied. We take $\mathrm{F}$ as the closure of $\BR$-rational points of
$\Spec (A)$ when $A$ satisfies (ii) and empty set otherwise. Let $\ma\subset A$ be an ideal such that $\dim(A/\ma)\le 1$ and $\dim(\ma \cap \mathrm{F})\le d-2$. Choose $\MI$ such that $\mathrm{F}=V(\MI)$.

 Now if $A$ satisfies (i) then $A/\ma$ is a real affine algebra of dimension $\le d-1$ such that $A/\ma$ does not have any real maximal ideal. In this case, it follows from Theorem \ref{PC} that all projective $A/\ma$-modules of rank $\ge d-1$ are cancellative. In this case $\MI=A$. Hence $\dim(V(J)\cap F)=\dim(V(A))\le d-2.$ The other hypotheses of \cite[Theorem 1.3]{M} follow from (a), (b), (c), and (III).
 
    Now we assume that $A$ satisfies (ii). Therefore, we have $\hh(\MI)\ge 2$. In this case if $\hh(\ma)\ge 2$ then by Bass cancellation \cite{HBASS} all projective $A/\ma$-module of rank $\ge d-1$ are cancellative. So we assume that $\hh(\ma)=1$. Since $A$ satisfies (ii) it follows that $A/\ma$ is a real affine algebra of dimension $d$ such that either $A/\ma$ does not have any real maximal ideals or the intersection of all real maximal ideals in $A/\ma$ has height at least $1$. Therefore, again applying Theorem \ref{PC} it follows that that all projective $A/\ma$-modules of rank $\ge d-1$ are cancellative. As $\hh(\MI)\ge 2$, we get $\dim(V(J)\cap F)=\dim(V(J+\MI))\le d-2.$ The other hypotheses of \cite[Theorem 1.3]{M} follow from (a), (b), (c), and (III).
    
    Therefore, applying \cite[Theorem 1.3]{M} there exist (A) a projective $A$-module $P$ of rank $d$ such that $J$ is a surjective image of $P$, and (B) $P\oplus Q\cong A^{2d}$. Now from (B), we get the following.
    $$[P]-[ A^{d}]=[A^{d}]-[Q]=-[A/K]\in \mathrm{K}_0(R)$$
    This concludes the proof.\qed

% 	\bd
%d-th Chern class $C_d(P)$
% \ed
% 
% \bc
% Let $R$ be reduced and let $P$ be a projective $R$-module of rank $d$. Then $P$ has a unimodular element if and only if d-th Chern class $C_d(P)$ is $0$ in $\text{CH}_0(\Spec(R))$. 
% \ec
% \proof Applying Swan's version of Bertini theorem [Theorem \ref{SwB}] we obtain a surjection $\alpha:P\surj I$, where $I\subset R$ is a reduced ideal of height $d$. Then we get $e(P)=(I)$ in $\E_0^d(R,L)$, where $\wedge^dP=L$. By Theorem \ref{QBRS} it follows that $C_d(P)=0$ in $\text{CH}_0(\Spec(R))$ if and only if $e(P)=0$ in $\E^d_0(R,L)$. Using \cite[Theorem 6.8]{SMBB3} we get that $e(P)=0$ in $\E^d_0(R,L)$ if and only if $(I)=0$ in $\E^d_0(R)$. It follows from [Proposition \ref{decg}, Step - 1, 2] and \cite[Corollary 4.4]{SMBB3} that $(I)=0$ in $\E^d_0(R,L)$ if and only if $P$ has a unimodular element. This concludes the proof.\qed
% 
 
 Let $A$ be a ring, and let $M$ be an $A$-module. We fix the following notations for the remaining part of this section.
                                                                                                                   
\begin{enumerate}[\quad\quad $\bullet$]
	
	\item  $supp(M)=\{p\in \Spec(A):M_p\not= 0\}$
	\item $\nu(p,M)=\mu(M_p)+\dim(A/p)$
	\item $ \eta(M)=\sup\{\nu(p,M):p\in supp(M)\}$
	\item $\delta(M)=\sup\{\nu(p,M):p\in supp(M),\,\dim(A/p)<\dim(A)\}$
	%   \item For any $m\in M$, $X(m)=\{p\in \Spec(A):m$ is a basic element of $M$ at $ p\}$;
	%   \item $\eta(M,m)=\sup\{\nu(p,M):p\in X(m)\}.$
\end{enumerate}

Now we state a theorem which is due to M. P. Murthy. We comment that with the results developed in this article so far, the exact same proof of \cite[Theorem 4.1]{M} will work in our set-up as well. To avoid repeating the same arguments we omit the proof.

\bp\label{MASL}
Let $A$ be as in Proposition \ref{PGM}, and let $M$ be an $A$-module. Then there is a projective $A$-module $P$ (with trivial determinant) of rank $\delta(M)$ such that (1) the module $M$ is a surjective image of $P$, and (2) $[P]-[A^{\delta(M)}]\in \mathrm{F^dK}_0(A)$.
\ep
   
   We now state and prove the main theorem of this section.
   
   \bt\label{EETT}
   Let $A$ be as in Proposition \ref{PGM}. Then the following are equivalent.
   \begin{enumerate}[\quad \quad (1)]
   	\item Every $A$-module $M$ is generated by at most $\delta(M)$ element.
   	\item Every locally complete intersection ideal in $A$ of height $d$ is complete intersection.
   	\item Every smooth maximal ideal in $A$ of height $d$ is complete intersection.
   	\item Every projective $A$-module of rank $d$ splits off a free factor of rank one.
   \end{enumerate}
   \et
   
   \proof $(1)\implies (2)$ If $I\subset A$ is locally complete intersection ideal in $A$ of height $d$, then $\delta(I)=d$. Hence, (2) follows.
   
    $(2)\implies (3)$ Since every smooth maximal ideal in $A$ of height $d$ is a locally complete intersection ideal of height $d$. Hence, (3) follows.
    
    $(3)\implies (4)$ Let $L$ be a projective $A$-module of rank one. Statement (3) will imply that $\widetilde{\E}^d_0(A)=0$. Combining Lemmas \ref{Step1}, \ref{Step2}, \ref{Step3} and \cite[Theorem 6.8]{SMBB3} we have the following isomorphisms:$$\widetilde{\E}^d_0(A)\cong {\E}^d_0(A)\cong {\E}^d_0(A,L)\cong {\E}^d(A,L).$$Hence it follows that ${\E}^d(A,L)=0$. Therefore, applying \cite[Corollary 4.4]{SMBB3} statement (4) follows.
    
    $(4)\implies (1)$ Let $I\subset A$ be an ideal such that $\mu(I/I^2)=\hh(I)=d$. Applying [Proposition \ref{PG2} (2)] there exists a projective $A$-module $P$ of rank $d$ such that $I$ is a surjective image of $P$. Since $P$ has a unimodular element, applying Corollary \ref{NMKL} we obtain that $\mu(I)=d$. This in particular proves that $\E^d_0(A)=0$. Hence using Corollary \ref{Cor}, we get that $\E^d(A)$ is trivial. Therefore, by Theorem \ref{k_0} we obtain that $\mathrm{F^dK_0}(A)=0$. Applying Proposition \ref{MASL} one can find a projective $A$-module $Q$ of rank $\delta(M)$ such that (a) the module $M$ is a surjective image of $Q$, and (b) $[Q]-[A^{\delta(M)}]\in \mathrm{F^dK}_0(A)$. If $\eta(M)\le \delta(M)$ then $M$ is generated by $\delta(M)$ elements \cite{Fo}. Therefore, we assume that $\eta(M)> \delta(M)$. In this case, one has $\delta(M)\ge d$, for details we refer to \cite[Proof of Theorem 4.1, paragraph 1]{M}. As $\mathrm{F^dK_0}(A)=0$ we get that $Q$ is stably free. Hence, by Theorem \ref{P} the module $Q$ is free. This concludes the proof. \qed
    
    We conclude this article with an example of an affine algebra over a field that satisfies the hypotheses of Theorem \ref{EETT}.
    
    \smallskip
    
    \example It follows from \cite[Theorem 4.4]{M} and \cite{AKMU}, that Theorem \ref{EETT} also holds for any affine algebra of dimension $d \ge 2$ over an algebraically closed field. Let $\overline{\mathbb{F}}_p$ denote the algebraic closure of the finite field $\mathbb{F}_p$, where $p$ is a prime. Let $A$ be an affine algebra of dimension $d \ge 2$ over $\overline{\mathbb{F}}_p$. Then, it is shown in \cite[Theorem 3.2]{MM} and \cite[Theorem 2.4]{MKD} that $\mathrm{E}^d(A) = 0$.

\bibliographystyle{abbrvurl}
%\bibliography{Realbib}

\end{document}